%% file: TQFTpost.tex
\documentclass[12pt]{article}

\title{Curves in Calabi-Yau threefolds and
Topological Quantum Field Theory} 
\author{Jim Bryan and Rahul Pandharipande} 
\date{\today}

\usepackage[PostScript=Rokicki]{diagrams}
\usepackage{epic}
\usepackage{eepic}

\usepackage{amsmath}
\usepackage{amsmath,amsthm,amsfonts}
\usepackage{verbatim}

\newtheorem{thm}{Theorem}[section]
\newtheorem{theorem}[thm]{Theorem}
\newtheorem{cor}[thm]{Corollary}

\newtheorem{lemma}[thm]{Lemma}

\newtheorem{proposition}[thm]{Proposition}

\newtheorem{definition}[thm]{Definition}
\newtheorem{warning}[thm]{Warning}

\newtheorem{rem1}[thm]{Remark}
\newenvironment{remark}{\begin{rem1}\em}{\end{rem1}}

\newcommand{\cnums} {{\mathbb C}}          
\newcommand{\qnums} {{\mathbb Q}}		
\newcommand{\affineline}{{\mathbb{A}^{1}}}                

\renewcommand{\P}{\mathbb{P}}
\newcommand{\M}{\overline{M}^{\bullet }}
\newcommand{\Mh}{\overline{M}_{h}}

\newcommand{\bfalpha}{{\boldsymbol{\alpha }}}
\newcommand{\bfbeta}{{\boldsymbol{\beta}}}
\newcommand{\bfgamma}{{\boldsymbol{\gamma  }}}
\newcommand{\bfdelta}{{\boldsymbol{\delta  }}}
\newcommand{\bfeta}{{\boldsymbol{\eta }}}

\newcommand{\W}{\mathfrak{W}}
\renewcommand{\O}{\mathcal{O}}
\newcommand{\X}{\mathcal{X}}

\newcommand{\combinatfactor}{\mathfrak{z}}
\newcommand{\Q}{{\mathbb Q}}
\newcommand{\com}{{\mathbb C}}

\newcommand{\one}{\mathbf{1}}
\newcommand{\Id}{\mathrm{Id}}

\begin{document}

\maketitle 

\begin{abstract}

We continue our study of the local Gromov-Witten invariants
of curves in Calabi-Yau threefolds.

We define relative invariants for the local theory which
give rise to a 1+1-dimensional TQFT taking values in the ring $\qnums
[[t]]$. The associated Frobenius algebra over $\qnums [[t]]$ is
semisimple. Consequently, we obtain a structure result for the local
invariants. As an easy consequence of our structure formula, we recover the
closed formulas for the local invariants in case either the target genus
or the degree equals 1.

\end{abstract}


\section{Notation, definitions and results}\label{sec: notation and results}

A central problem in Gromov-Witten theory is to determine the structure of
the Gromov-Witten invariants. Of special interest is the case where the
target manifold is a Calabi-Yau threefold. We prove a structure result for
the local Gromov-Witten invariants of a curve in a Calabi-Yau threefold.

First, we define a relative version of the local invariants. In
Theorem~\ref{thm: the TQFT exists}, we prove the new invariants determine a
1+1-dimensional semisimple TQFT taking values in the ring $\Q[[t]]$.  The
structure formula for the invariants, Theorem~\ref{thm: Zd(g) formula}, is
obtained from the semisimple TQFT and completely determines the dependence
of the invariants on the target genus.

The local invariants are derived from considering rigid curves in
Calabi-Yau threefolds. Let $X$ be a non-singular curve of genus $g$ in a
Calabi-Yau threefold $Y$. Assuming certain rigidity conditions on $X\subset
Y$, there are well defined local Gromov-Witten invariants of $X$ in $Y$.
That is, the contribution to the Gromov-Witten invariants of $Y$ by maps
with image $X$ is well defined. These local invariants depend only on $g$,
$h$, and $d$, which are respectively, the target genus, the domain genus,
and the degree of the maps to $X$.

In \cite{Br-Pa}, we defined an integral depending only on $d$, $h$, and $g$
that gives the value of the local invariants whenever $X\subset Y$
satisfies the requisite rigidity. The integral is given by
\[
N_{d}^{h-g} (g)=\int  _{[\Mh (X,d[X])]^{vir}} c (I (X)),
\]
see \cite{Br-Pa} (c.f. \cite{Br-Pa-rigidity}) for details of this
discussion.  Here, $c (\cdot )$ denotes total Chern class and $I (X)$ is
defined by

\begin{align*}
I (X)&= -R^{\bullet  }\pi _{*}f^{*} (K_{X}\oplus \O_{X})\\
&=\hspace{7pt} R^{1 }\pi _{*}f^{*} (K_{X}\oplus \O_{X})-R^{0 }\pi
_{*}f^{*} (K_{X}\oplus \O_{X})
\end{align*}
which is understood to be an element of $K$-theory.

It will be convenient to work with the version of the local invariants
corresponding to stable maps with possibly disconnected domains:

\begin{definition}\label{def: possibly disconnected moduli space}
Let $X$ be a non-singular curve of genus $g$. Let $\M (X,d)$ be the moduli
space of degree $d$ stable maps
$$f:C\to X$$ 
where $C$ is a possibly disconnected curve. 
We require $f$ to be nonconstant  on each
connected component of $C$. 
Let $$[\M (X,d)]^{vir}\in
A_{*} (\M (X,d);\qnums )$$ denote the virtual fundamental cycle of
the moduli space.  
\end{definition}

Following \cite{Okounkov-Pandharipande-completed-cycles}, the superscript
$^{\bullet }$ is used to denote the moduli space with possibly disconnected
domain curves.  The usual genus subscript is omitted as we consider all
domain genera: the moduli space is a countable union of connected
components with varying expected dimensions.  The branch points of a stable
map to $X$ are well-defined by \cite{Fantechi-Pandharipande}.  The number
of branch points $b$ of map $f$ equals the expected dimension of the moduli
space at the moduli point $[f]$.

We define the \textbf{(possibly disconnected) local Gromov-Witten
invariants} to be
\[
Z_{d}^{b} (g)=\int _{[\M (X,d)]^{vir}} c_{b} (I (X)).
\]

The relationship between the possibly disconnected invariants and the
(connected) invariants $N^{h-g}_{d} (g)$ is easily seen to be

\[
\sum _{d>0}\sum _{b\geq 0}Z^{b}_{d} (g)t^{b}q^{d} =
\operatorname{exp}\left( \sum _{d>0}\sum _{b\geq 0}N^{h-g}_{d}
(g)t^{b}q^{d}\right)
\]
where
$$2h-2= (2g-2)d +b .$$
The generating function for the degree $d$, local, disconnected invariants is
\[
Z_{d} (g)=\sum _{b=0}^{\infty }Z_{d}^{b} (g)t^{b}.
\]

The series $Z_{d} (g)$ is our basic object of study. Clearly, the
disconnected invariants $Z_{d}^{b} (g)$ and the connected invariants
$N^{h-g}_{d} (g)$ contain equivalent information.

In Section \ref{sec: relative local invs}, we use J. Li's theory of
relative stable maps \cite{Li-relative1,Li-relative2} to construct relative
versions of the local invariants. These relative invariants obey a gluing
law which allows us to construct a Topological Quantum Field Theory (TQFT):

\begin{theorem}\label{thm: the TQFT exists}
There exists a $1+1$-dimensional TQFT, $Z_d(-)$, with the following three properties:
\begin{enumerate}
\item[(i)] $Z_d(-)$ is semisimple,
\item[(ii)] $Z_d(-)$ takes values in $\Q[[t]]$,
\item[(iii)] $Z_d(-)$ applied to a genus $g$ closed surface yields the value $Z_d(g)$,
             the generating series of the local invariants.
\end{enumerate}
\end{theorem}

The $t=0$ specialization of $Z_d(-)$ 
is a well-known TQFT obtained from the gauge theory
of the symmetric group $S_{d}$, see Lemma~\ref{lemma: 
Z at t=0 is the
DW/FQ TQFT}. The TQFT determined by $S_d$ 
was studied by Dijkgraaf-Witten and Freed-Quinn
\cite{Dijkgraaf-Witten90,Freed-Quinn}. Our TQFT may be viewed as a 
$1$-parameter deformation of the Dijkgraaf-Witten/Freed-Quinn theory.

Corresponding to any 1+1-dimensional TQFT is a Frobenius algebra. In our
case, the dimension of the corresponding Frobenius algebra is $p
(d)$, the number of partitions of $d$. As a corollary of Theorem 
\ref{thm: the TQFT exists}, we
deduce the following structure formula.

\begin{theorem}\label{thm: Zd(g) formula}
There exist universal power series $\lambda _{\alpha }\in \qnums [[t]]$,
labelled by partitions $\alpha $ of $d$ (denoted $\alpha \vdash d$), which
determine the local invariants by:
\[
Z_{d} (g)=\sum _{\alpha \vdash d} \lambda _{\alpha }^{g-1}.
\]
Moreover, the constant term of the series $\lambda _{\alpha }$ is given by
\[
\left(\frac{d!}{\dim R_{\alpha }} \right)^{2}
\]
where $R_{\alpha }$ is the irreducible representation of the symmetric
group associated to $\alpha $.
\end{theorem}

The two main theorems of \cite{Pandharipande-degenerate-contributions}
(theorems 1 and 2) compute the local invariants in the case of degree 1 and
the case of target genus 1. We recover these two results as immediate
corollaries of the above structure theorem (Corollaries \ref{cor: d=1
follows from structure thm} and \ref{cor: super-rigid elliptic
contributions} below):

\begin{cor}\label{cor: super-rigid elliptic contributions}
The series $Z_{d} (1)$ is the constant series $p (d)$. In particular, the
genus two and higher multiple cover contributions of a super-rigid elliptic
curve are all zero.
\end{cor}

By the localization calculation of Faber-Pandharipande \cite{Fa-Pa}, $Z_{d}
(0)$ is given by
\[
Z_{d} (0)= \sum _{\alpha \vdash d}\frac{t^{2d}}{\combinatfactor (\alpha )}
\prod _{i=1}^{\ell(\alpha) } \left(2\sin (\frac{\alpha _{i}t}{2}) \right)^{-2}
\]
where the sum is over all partitions $\alpha$ of $d$. Here,
$\ell(\alpha) $ is the length of the $\alpha$, and $\combinatfactor (\alpha )$ is a
combinatorial factor (see Definition~\ref{defn: e(a) where a is a
partition}). In particular, for $d=1$, we have
\[
Z_{1} (0)= \left(\frac{\sin (t/2)}{t/2} \right)^{-2}
\]
which we combine with our structure formula to deduce the following:
\begin{cor}\label{cor: d=1 follows from structure thm}
The $d=1$ local invariants are given by
\[
Z_{1} (g)=\left(\frac{\sin (t/2)}{t/2} \right)^{2g-2}.
\]
\end{cor}

Recent progress (to be explained in a forthcoming paper
\cite{Bryan-Pandharipande-in-prep}) has allowed us to completely determine
the TQFT $Z_{d} (-)$ for small values of $d$. In particular, for $d=2$ we have
\begin{theorem}[\cite{Bryan-Pandharipande-in-prep}]\label{thm: d=2 formula}
The $d=2$ local invariants are given by
\[
Z_{2} (g)=\left(\frac{\sin (t/2)}{t/2} \right)^{4g-4}\left\{(4-4\sin
(t/2))^{g-1}+ (4+4\sin (t/2))^{g-1} \right\}.
\]
\end{theorem}

The above formula gives the double cover contributions of any
smooth curve in a Calabi-Yau threefold with a generic normal bundle. The above
formula also verifies the local Gopakumar-Vafa conjecture for degree 2
maps, i.e. the corresponding BPS invariants are integers. See \cite{Br-Pa}
for a discussion of the BPS invariants and the local Gopakumar-Vafa
conjecture.

\section{Semisimple TQFTs over complete local rings}\label{sec: TQFTs}
Let $(n+1)\mathbf{Cob}$ be the symmetric monoidal category with
objects given by compact oriented $n$-manifolds and morphisms 
given by (diffeomorphism
classes of) oriented cobordisms.
An \textbf{$(n+1)$-dimensional TQFT with values in a commutative ring
$R$} is a symmetric monoidal functor
\[
Z:(n+1)\mathbf{Cob}\to R\mathbf{mod},
\]
where $R\mathbf{mod}$ is the category of $R$-modules.
The definition amounts to the following axioms for $Z$:
\begin{enumerate}
\item[(i)] To each compact oriented $n$-manifold $Y$, $Z$
assigns an $R$-module $Z(Y)$.
\item[(ii)]  To each oriented cobordism $W$ from $Y_{1}$ to $Y_{2}$, $Z$ assigns
an $R$-module homomorphism $Z (W):Z (Y_{1})\to Z (Y_{2})$.
\item[(iii)] If two oriented cobordisms are equivalent 
$W\cong W'$ by a boundary preserving diffeomorphism, then $Z (W)=Z (W')$.
\item[(iv)] The trivial oriented cobordism corresponds to the identity
homomorphism, $Z (Y\times
[0,1])=\Id_{Z (Y)}$.
\item[(v)] The concatenation of cobordisms corresponds to the composition of
the corresponding $R$-module homomorphisms.
\item[(vi)] The disjoint union of $n$-manifolds corresponds to the tensor product of
$R$-modules, $Z (Y_{1}\coprod Y_{2})=Z (Y_{1})\otimes Z (Y_{2})$, and the
disjoint union of cobordisms corresponds to the tensor product of homomorphisms, $Z
(W_{1}\coprod W_{2})=Z (W_{1})\otimes Z (W_{2})$.
\item[(vii)] The empty $n$-manifold corresponds to the ground ring, $Z (\emptyset
)=R$.
\end{enumerate}

A compact oriented $(n+1)$-manifold $W$ may be viewed as a oriented cobordism 
between empty manifolds. Then,
$$Z(W)\in Hom_R(R,R) \cong R.$$
The element $Z(W) \in R$ is a topological invariant of $W$. 


TQFTs of dimension 1+1 are in bijective correspondence with commutative
Frobenius algebras. The result goes back to Dijkgraaf's thesis, and has
been proven in various contexts by Sawin \cite{Sawin}, Abrams
\cite{Abrams}, and Quinn \cite{Quinn}. The form of the correspondence that we
quote is due to Kock \cite{Kock:FA-2DTQFT}:

\begin{theorem}\label{thm: 1+1 TQFTs = comm Frob algs}
The category of 1+1-dimensional TQFTs taking values in $R$ is
equivalent to the category of commutative Frobenius algebras over $R$.
\end{theorem}

A \textbf{commutative Frobenius algebra over $R$ } is a
commutative $R$-algebra $A$ equipped with a counit $\mu :A\to R$ and a
coassociative, cocommutative, comultiplication $\triangle :A\to A\otimes A$
satisfying the Frobenius relation and the counit axiom:
\begin{align*}
(m\otimes \Id) (a\otimes \triangle (b))&= (\Id\otimes m) (\triangle (a)\otimes b)=\triangle (m (a\otimes b))\\
(\Id\otimes \mu) (\triangle (a))&= (\mu \otimes \Id) (\triangle (a))=a
\end{align*}
where $m:A\otimes A\to A$ is multiplication. The axioms imply that $A$ is
finitely generated as an $R$-module.

Given an invertible element $\lambda \in R$, we can give $R$ the structure
of a Frobenius algebra by setting $\triangle (1)=\lambda $ and
(consequently) $\mu (1)=\lambda ^{-1}$. We denote this Frobenius algebra by
$R_{\lambda }$. A Frobenius algebra is \textbf{semisimple} if it is
isomorphic to $R_{\lambda _{1}}\oplus \dots \oplus R_{\lambda _{n}}$ for
some $\lambda _{1},\dots ,\lambda _{n}\in R$.

A 1+1-dimensional TQFT is semisimple if the corresponding Frobenius
algebra is semisimple. 
If $W_{g} $ is a closed surface of genus $g$, and $Z$ is a
semi-simple TQFT, then an elementary argument from the axioms yields:
\begin{equation}\label{eqn: formula for invariant in a semi-simple TQFT}
Z (W_{g})=\sum _{i=1}^{n}\lambda _{i}^{g-1}.
\end{equation}

The following basic result is the key to proving the semisimplicity of the
TQFT that we construct from local Gromov-Witten invariants.

\begin{proposition}\label{prop: A/m semi-simple over R/m => A s-s over R}
Let $R$ be a complete local ring. Let $m\subset R$ be the maximal ideal and
let $A$ be a Frobenius algebra over $R$. Suppose that $A$ is free as an
$R$-module and that $A/mA$ is a semi-simple Frobenius algebra over
$R/m$. Then $A$ is semi-simple (over $R$).
\end{proposition}

\textsc{Proof}: Let $e_{1},\dots ,e_{n}\in A$ be representatives for an
idempotent basis of $A/mA$, that is $e_{i}e_{i}-e_{i}\in m$ for all $i$ and
$e_{i}e_{j}\in m$ for all $i\neq j$. By Nakayama's lemma, $\{e_{i} \}$ is a
basis of $A$, and we wish to construct a new basis which is idempotent in
$A$. We begin by constructing an idempotent basis in $A/m^{2}A$.

It is easy to see that an element is invertible in $A$ if and only if it is
invertible in $A/m$. In particular, $1-2e_{i}$ is invertible since its
square is $1$ modulo $m$.  Let $b_{i}=e_{i}e_{i}-e_{i}$ and set
\[
e_{i}'=e_{i}+b_{i} (1-2e_{i})^{-1}.
\]
A short computation shows that
\[
e_{i}'e_{i}'-e_{i}'=b_{i}^{2} (1-2e_{i})^{-2}
\]
which is in $m^{2}$ for all $i$ since $b_{i}\in m$. Then $e'_{i}e'_{j}\in
m^{2}$ for $i\neq j$ follows from the facts that $(e'_{i}e'_{j})^{2}$,
$e'_{i}e'_{i} (e'_{j}e'_{j}-e'_{j})$, and $e'_{j} (e'_{i}e'_{i}-e'_{i})$
are all in $m^{2}$. Thus $\{e_{i}' \}$ is an idempotent basis for
$A/m^{2}$.

We construct a sequence of bases $\{e'_{i} \},\{e''_{i} \},\dots
,\{e_{i}^{(k)} \}$ by setting
\begin{align*}
b_{i}^{(k)}&=e_{i}^{(k)}e_{i}^{(k)}-e_{i}^{(k)}\\
e_{i}^{(k+1)}&=e_{i}^{(k)}+b_{i}^{(k)} (1-2e_{i}^{(k)})^{-1}.
\end{align*}
The same argument as above shows that $\{e_{i}^{(k)} \}$ is an idempotent
basis for $A/m^{k+1}A$. Since $R$ is complete, there exists
$\tilde{e}_{i}\in A$ such that $\tilde{e}_{i}=e_{i}^{(k)}\bmod m^{k+1}$ for
all $k$. By construction, $\tilde{e}_{i}$ is an idempotent basis for
$A$. Let $\lambda _{i}=\mu (\tilde{e}_{i})^{-1}$. Since $A$ is free as an
$R$-module, each $\tilde{e}_{i}$ generates an $R$ summand of $A $ and thus
we have constructed the desired isomorphism of Frobenius algebras:
\[
A\cong R_{\lambda _{1}}\oplus \dots \oplus R_{\lambda _{n}}.
\]
\qed

\section{Relative local invariants and gluing}\label{sec:
relative local invs}

Motivated by the symplectic theory of A.-M. Li and Y. Ruan \cite{Li-Ruan},
J. Li has developed an algebraic theory of relative stable maps to a pair
$(X,B)$. This theory compactifies the moduli space of maps to $X$ with
prescribed ramification over a non-singular divisor $B\subset X$,
\cite{Li-relative1,Li-relative2}.  Li constructs a moduli stack of relative
stable maps together with a virtual fundamental cycle and proves a gluing
formula.  Consider a degeneration of $X$ to $X_{1}\cup_{B} X_{2}$, the
union of $X_{1}$ and $X_{2}$ along a smooth divisor $B$. The gluing formula
expresses the virtual fundamental cycle of the usual stable map moduli
space of $X$ in terms of virtual cycles for relative stable maps of
$(X_{1},B)$ and $(X_{2},B)$. The theory of relative stable maps has also
been pursued in \cite{Ionel-Parker00, Ionel-Parker-Annals2003}, \cite{EGH}.

In our case, the target is a non-singular curve $X$ of genus $g$,
and the divisor $B$ is a collection of points $b_{1},\dots ,b_{r}\in X$.

\begin{definition}\label{defn: rel stable maps}
Let $(X,b_{1},\dots b_{r})$ be a fixed non-singular genus $g$ curve with
$r$ distinct marked points.  Let $\alpha_1, \ldots, \alpha_r$ be
partitions of $d$.
Let $$\M (X, (\alpha _{1}\dots \alpha_r))$$ be 
the moduli stack of relative stable maps (in the sense of
 Li)\footnote{For a formal definition of relative stable maps, we refer to
\cite{Li-relative1} Section~4.} with target $(X,b_{1},\dots ,b_{r})$ satisfying the
following:
\begin{enumerate}
\item[(i)] The maps have degree $d$.
\item[(ii)] The maps are ramified over $b_{i}$ with ramification type $\alpha
_{i}$.
\item[(iii)] The domain curves are possibly disconnected, but the map is not
degree 0 on any connected component.
\item[(iv)] The domain curves are \emph{not} marked.
\end{enumerate}
The partition $\alpha _{i}\vdash d$ determines a ramification type over
$b_{i}$ by requiring the monodromy of the cover (considered as a conjugacy
class of $S_{d}$) has cycle type $\alpha _{i}$.  We will use a multi-index
$\bfalpha = (\alpha _{1},\dots ,\alpha _{r})$ to shorten the notation to
$\M (X, {\bfalpha})$ when there is no risk of confusion.
\end{definition}

Our moduli spaces of relative stable maps 
differ from  Li's in a few minor ways:
\begin{enumerate}
\item[(i)] We do not require the domain to be connected.
\item[(ii)] We do not fix the domain genus, so $\M (X ,\bfalpha)$
is a countable union of components.
\item[(iii)] We do not mark the domain curve at all. Li's spaces include a
marking of the ramification locus.
\end{enumerate}
We make these modifications to Li's theory to simplify
 the combinatorics of
the gluing theory. It is straightforward to express our moduli
spaces in terms of unions, products, and finite quotients of Li's
spaces.

There is a universal diagram (of stacks):
\[
\begin{diagram}
R & \rInto & U & \rTo^{f} & \X & \lInto &B\\
&&\dTo>{\pi }&\ldTo>{p}&&&\\
&&\M (X,\bfalpha  )&&&&
\end{diagram}
\]
where $U$ is the universal domain curve, $\X$ is the universal target
curve, $f$ is the universal map, $B$ is the universal prescribed branch
divisor, and $R$ is the universal prescribed ramification divisor. The
divisors $B$ and $R$ are taken with reduced structure. Note that the map
$f$ can also be ramified away from $R$, but at $R$ the map $f$ has the
ramification type prescribed by the data $\bfalpha $.

One of the salient features of the relative theory is that the target curve
may ``bubble'' off rational components meeting the original $X$ in
nodes. That is to say, the family
$$p:\X\to \M (X, \bfalpha)$$ is nontrivial: special fibers have chains of
rational curves attached to $X$ at the points $b_{i}$. However, the
universal prescribed branch locus $B$ lies in the non-singular locus of the
fibers of $p$.

Let $[\M (X ,\bfalpha )]^{vir}$ denote the virtual fundamental class. With
our conventions, the virtual class is a countable sum of cycle classes with
degree $b$ part supported on the components of expected dimension $b$. The
expected dimension of a component is given by
\[
b= -\chi + \ell(\bfalpha) -d (2g-2+r)
\]
where $\ell(\bfalpha) =\ell(\alpha _{1})+\dots +\ell(\alpha _{r})$ 
is the sum of the
lengths of the partitions and $\chi $ is the domain Euler characteristic.

Let $\omega $ be the relative dualizing sheaf of $p$. 
We define an element $I (X,\bfalpha )$ of the $K$-theory
of coherent sheaves on $\M (X, \bfalpha )$ by
\begin{align*}
I (X,\bfalpha )&=-R^{\bullet }\pi _{*} (f^{*} (\omega
(B))\oplus
\O (-R))\\
&=\hspace{7pt} R^{1 }\pi _{*} (f^{*} (\omega (B))\oplus
\O (-R))-R^{0 }\pi _{*} (f^{*} (\omega
(B))\oplus \O (-R)).
\end{align*}

We define the \textbf{relative local invariants} by
\[
Z^{b}_{d} (g)_\bfalpha=\int _{[\M (X,\bfalpha  )]^{vir}}c_{b}
(I (X, \bfalpha ))
\]
and their corresponding generating series $Z_{d} (g)_\bfalpha \in \qnums
[[t]]$ by
\[
Z_{d} (g)_\bfalpha =\sum _{b=0}^{\infty }Z_{d}^{b} (g)_\bfalpha \; t^{b}.
\]
The element $I (X, \bfalpha )$ has rank $b$ when restricted to the
components of the moduli space $\M (X, \bfalpha )$ with expected dimension $b$. We
regard $I (X,\bfalpha )$ as a (virtual) obstruction bundle and $c_{b}
(I (X, \bfalpha) )$ as the associated Euler class.

\begin{remark}\label{rem: open string interpretation}
In case $g=0$ and $r=1$, the above Euler class was previously defined by
J. Li and Y. Song \cite{Li-Song}. Li and Song were modeling ``open string''
Gromov-Witten theory using relative stable maps. They argued that the
integral $Z_{d}^{b} (0)_{\alpha }$ should compute the multiple cover
formula for maps to a disk in a Calabi-Yau threefold where the boundary of the
disk lies on a fixed Lagrangian 3-manifold (see also Katz-Liu
\cite{Katz-Liu}). They arrived at the integrand by considering the
obstruction theory of maps to the disk with the Lagrangian boundary
conditions. It is plausible that such an open-string theory interpretation
of the relative local invariants exists in general.
\end{remark}

We will use Li's theory to obtain  ``gluing relations'' among the
relative (and usual) local invariants.
The following combinatorial quantity arises frequently:

\begin{definition}\label{defn: e(a) where a is a partition}
Let $\gamma \vdash d$ be a partition of $d$ and let $\gamma (k)$ be the
number of parts of size $k$ in $\gamma$, so  $d=\sum _{k=1}^{\infty
}\gamma (k)k$. We define:
\[
\combinatfactor (\gamma )=\prod _{k=1}^{\infty }k^{\gamma (k)}\gamma (k)!.
\]
If $c (\gamma )\subset S_{d}$ denotes the conjugacy class in
the symmetric group consisting of elements having cycle type $\gamma $, then
$\combinatfactor (\gamma )$ is the order of the centralizer of $c (\gamma )$.
\end{definition}

The basic gluing laws are given by the following:

\begin{theorem}\label{thm: gluing laws}
Let $\bfalpha= (\alpha_1, \ldots, \alpha_r)$.
For any choice $g_{1}+g_{2}=g$ and any splitting $$\{\alpha _{1},\dots
,\alpha _{r} \}=\{\alpha _{1},\dots ,\alpha _{k} \}\cup \{\alpha
_{k+1},\dots ,\alpha _{r} \},$$ we have:
\begin{equation}\label{eqn: gluing along separating curve}
Z_{d} (g)_\bfalpha=
\sum _{\gamma \vdash d}\combinatfactor
(\gamma )Z_{d} (g_{1})_{\alpha _{1},\dots, \alpha _{k},\gamma }\; Z_{d}
(g_{2})_ {\alpha _{k+1},\dots, \alpha _{r},\gamma}.
\end{equation}
We also have
\begin{equation}\label{eqn: gluing along non-separating curve}
Z_{d} (g+1)_\bfalpha=\sum _{\gamma \vdash d}\combinatfactor
(\gamma )Z_{d} (g)_{\alpha _{1},\dots, \alpha _{r},\gamma, \gamma }.
\end{equation}
\end{theorem}

The first formula corresponds to splitting a genus $g$ surface with
$r$ boundaries along a separating curve to obtain two surfaces of genus
$g_{1} $ and $g_{2}$ with $(k+1)$ and $(r-k+1)$ boundaries. The second
formula corresponds to cutting a genus $g+1$ surface with $r$ boundaries
along a non-separating curve to obtain a genus $g$ surface with $(r+2)$
boundaries. We defer the proofs of these formulas to Appendix~\ref{appendix: proof of
gluing}.

\section{The TQFT}\label{sec: the TQFT}

We will show the gluing formulas of Theorem~\ref{thm:
gluing laws} allow us to organize the invariants $Z_{d} (g)_\bfalpha $
into a 1+1-dimensional TQFT over $\qnums [[t]]$. Throughout Section
\ref{sec: the TQFT},
$R$ will denote the ring of formal power series in $t$ over $\qnums $:
\[
R=\qnums [[t]].
\]

To define our TQFT, $Z_{d} (-)$, we need an $R$-module $H=Z_{d} (S^{1})$
associated to the circle ($H$ is the ``Hilbert space'' of the theory). We
define
\[
Z_{d} (S^{1})=H=\bigoplus  _{\alpha \vdash d}R\ e_{\alpha }
\]
to be a free $R$-module with a basis $\{e_{\alpha } \}_{\alpha \vdash d}$
labeled by partitions of $d$.

Using the given basis for $H$, we can express any module homomorphism 
\[
f:H^{\otimes r}\to H^{\otimes s}
\]
in tensor notation $f_{\alpha _{1}\dots \alpha _{r}}^{\beta
_{1}\dots \beta _{s}}\in R$ by
\[
f(e_{\alpha _{1}}\otimes \dots \otimes e_{\alpha _{r}}) = \sum _{\beta
_{1},\dots ,\beta _{s}}f_{\alpha _{1}\dots \alpha _{r}}^{\beta _{1}\dots
\beta _{s}}\ e_{\beta _{1}}\otimes \dots \otimes e_{ \beta _{s}}.
\]
Using multi-index notation and the Einstein summation convention,
we simply write:
\[
f:e_{\bfalpha} \mapsto f^{\bfbeta }_{\bfalpha }e_{\bfbeta} 
\]

We raise indices by the following formula: 
\[
Z_{d} (g)^{\beta _{1},\dots,\beta _{s}}_{\alpha _{1},\dots, \alpha _{r}} =
\combinatfactor (\beta _{1})\cdots \combinatfactor (\beta _{s})Z_{d} (g)_{\alpha
_{1},\dots, \alpha _{r},\beta _{1},\dots, \beta _{s}}.
\]
Then the gluing laws can be written succinctly:
\begin{align*}
Z_{d} (g_{1}+g_{2})^{\bfbeta, \bfdelta }_{\bfalpha ,\bfeta}&= Z_{d}
(g_{1})^{\bfbeta, \gamma }_{\bfalpha }Z_{d} (g_{2})^{\bfdelta }_{\bfeta,
\gamma }\\
Z_{d} (g+1)^{\bfbeta }_{\bfalpha }&=Z_{d} (g)_{\bfalpha ,\gamma }^{\bfbeta,
\gamma }
\end{align*}

Note that in the above equation, $\gamma $ it is a single index, while the
boldface indices ($\bfalpha $, $\bfbeta $, etc) stand for
multi-indices. Since the $\gamma $ is repeated on the left hand side, it is
summed over by convention.

Let $W_{r}^{s} (g)$ be the connected, oriented, genus g, cobordism from a
disjoint union of $r$ boundary circles to $s$ boundary circles. We define
\[
Z_{d} (W^{s}_{r} (g)):H^{\otimes r}\to H^{\otimes s}
\]
by
\[
e_{\bfalpha} \mapsto Z_{d} (g)_{\bfalpha }^{\bfbeta }e_{\bfbeta }
\]
where $\bfalpha =\alpha _{1}, \dots, \alpha _{r}$ and $\bfbeta =\beta _{1},
\dots , \beta _{s}$. 

For a disconnected cobordism $W=W_{1}\sqcup \dots
\sqcup W_{n}$, we define
\[
Z_{d} (W)=Z_{d} (W_{1})\otimes \dots \otimes Z_{d} (W_{n}).
\]

\begin{proposition}\label{prop: gluing formulas imply Z is a TQFT}
The functor $Z_{d} (-)$ defined above is a $(1+1)$-dimensional TQFT over
$R$.
\end{proposition}

\textsc{Proof:} To show $Z_{d} (-)$ is indeed a functor, we must 
prove that $Z_{d} (-)$ takes the concatenation of cobordisms to the
composition of $R$-module homomorphisms.
The composition 
\[
Z_{d} (W^{t}_{s} (g_{2}))\circ Z_{d} (W_{r}^{s} (g_{1})):H^{\otimes r}\to
H^{\otimes s}\to H^{\otimes t}
\]
determined by connected cobordisms is given by
\[
e_{\bfalpha} \mapsto Z_{d} (g_{1})^{\bfbeta }_{\bfalpha }e_{\bfbeta}\mapsto
Z_{d} (g_{1})^{\bfbeta }_{\bfalpha }Z_{d} (g_{2})^{\bfgamma }_{\bfbeta
}e_{\bfgamma}
\]
for $$\bfalpha =(\alpha _{1}, \dots, \alpha _{r}), \ \ 
\bfbeta =(\beta _{1},
\dots, \beta _{s}), \  \  \bfgamma =(\gamma _{1}, \dots , \gamma
_{t}).$$ 
Applying the gluing laws we obtain
\begin{align}\label{eqn: gluing for general connected cobod}
Z_{d} (g_{1})_{\bfalpha }^{\beta _{1}\dots \beta _{s}}Z_{d}
(g_{2})^{\bfgamma }_{\beta _{1}\dots \beta _{s}}&= Z_{d}
(g_{1}+g_{2})^{\bfgamma \beta _{2}\dots \beta _{s}}_{\bfalpha \beta
_{2}\dots \beta _{s}}\nonumber \\
&=Z_{d} (g_{1}+g_{2}+1)^{\bfgamma \beta _{3}\dots \beta _{s}}_{\bfalpha
\beta _{3}\dots \beta _{s}}\nonumber \\
&\quad \vdots \nonumber \\
&=Z_{d} (g_{1}+g_{2}+s-1)^{\bfgamma }_{\bfalpha } \ .
\end{align}
We have proven 
\[
Z_{d} (W^{t}_{s} (g_{2}))\circ Z_{d} (W_{r}^{s} (g_{1}))=Z_{d} (W^{t}_{r}
(g_{1}+g_{2}+s-1)).
\]
Since the concatenation of $W^{s}_{r} (g_{1})$ followed by $W_{s}^{t}
(g_{2})$ is $W^{t}_{r} (g_{1}+g_{2}+s-1)$ (see the Figure),
we have shown that $Z_{d} (-)$ is a functor, at least when applied to the
subcategory of $2\mathbf{Cob}$ consisting of connected cobordisms.

\begin{figure}\label{fig: cobord}
\centering 
\input{cob4.eepic} 
\caption{$W^{s}_{r} (g_{1})$ concatenated
with $W_{s}^{t} (g_{2})$ is $W^{t}_{r} (g_{1}+g_{2}+s-1)$.  The gluing
formula expressed by this picture is given by Equation ~(\ref{eqn: gluing
for general connected cobod}).}
\end{figure}
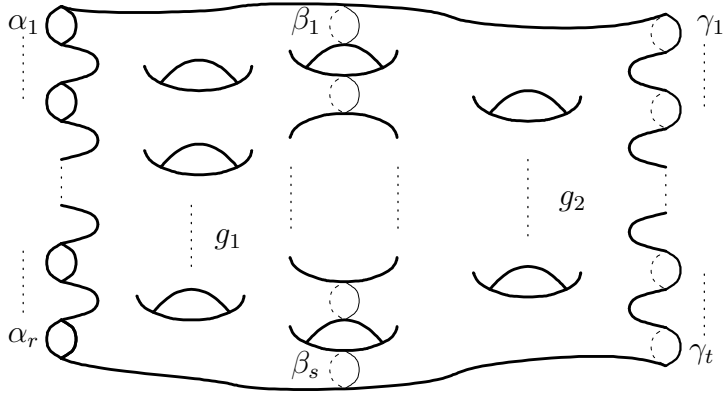

Similar computations apply to concatenations of disconnected
cobordisms. For example, the concatenation of the cobordism $W_{r}^{1}
(g_{1})\sqcup W^{1}_{s} (g_{2})$ followed by $W_{2}^{t} (g_{3})$ yields the
cobordism $W_{r+s}^{t} (g_{1}+g_{2}+g_{3})$. Correspondingly, for 
$$\bfalpha
=(\alpha _{1}, \dots , \alpha _{r}), \ \ \bfbeta =(\beta _{1}, \dots , \beta
_{s}), \ \ \bfdelta =(\delta _{1}, \dots , \delta _{t}),$$ the composition
\[
Z_{d} (W^{t}_{2} (g_{3}))\circ Z_{d} (W_{r}^{1} (g_{1})\sqcup W_{s}^{1}
(g_{2}))
\]
is given by 
\begin{align*}
e_{\bfalpha} \otimes e_{\bfbeta} &\mapsto Z_{d} (g_{1})^{\gamma
_{1}}_{\bfalpha }Z_{d} (g_{2})^{\gamma _{2}}_{\bfbeta }e_{\gamma
_{1}}\otimes e_{ \gamma _{2}}\\
&\mapsto Z_{d} (g_{1})^{\gamma _{1}}_{\bfalpha }Z_{d} (g_{2})^{\gamma _{2}}_{\bfbeta }Z_{d} (g_{3})^{\bfdelta }_{\gamma _{1}\gamma _{2}}e_{\bfdelta} \\
&=Z_{d} (g_{1})^{\gamma _{1}}_{\bfalpha }Z_{d} (g_{2}+g_{3})^{\bfdelta
}_{\bfbeta \gamma _{1}} e_{\bfdelta}\\
&=Z_{d} (g_{1}+g_{2}+g_{3})^{\bfdelta }_{\bfalpha \bfbeta }e_{\bfdelta}
\end{align*}
The general concatenation of cobordisms can be checked with similar
computations.

To prove that $Z_{d} (-)$ is a symmetric monoidal functor, we must
also
check that $Z_{d} (-)$ takes the trivial cobordism $S^{1}\times [0,1]$
to the identity. This is equivalent to the following lemma.

\begin{lemma}\label{lem: local invs of the annulus}
The local invariant of $\P ^{1}$, relative 2 points is given by:
\[
Z_{d} (0)_{\alpha \beta }=
\begin{cases}
\frac{1}{\combinatfactor (\alpha )}&\text{ if }\alpha =\beta, \\
0& \text{ if } \alpha \neq \beta .
\end{cases}
\]
\end{lemma}

\textsc{Proof:} The relative invariant $Z_{d} (0)_{\alpha \beta} $ can be
computed by virtual localization since there is a $\cnums ^{\times }$
action on $\P ^{1}$ preserving the relative points $0,\infty \in \P^1$ and
hence $\cnums ^{\times }$ acts on $\M (\P ^{1},(\alpha ,\beta ))$. However,
we can compute $Z_{d} (0)_{\alpha \beta }$ more easily as follows.

The component of $\M (\P ^{1},( \alpha, \beta ))$ of
virtual dimension 0 parameterizes stable maps $f:C\to \P ^{1}$ that are
unramified away from the prescribed ramification points $0,\infty \in \P
^{1}$. Any such map must be of the form
\[
f:\P ^{1}\sqcup \dots \sqcup \P ^{1}\to \P ^{1}
\]
where on the $i$th component, $f$ is of the form $z\mapsto z^{\alpha _{i}}$
for some $\alpha \vdash d$. 

Therefore, if $\alpha \neq \beta $, then the virtual dimension 0
component of the moduli space 
$\M (\P ^{1}, (\alpha, \beta))$ is empty. If $\alpha =\beta
$, then virtual dimension 0 component
 consists of a single moduli point $[f]$ corresponding to the
above map. The map $f$ has an automorphism group of order $\combinatfactor (\alpha
)$, hence
\[
Z^{0}_{d} (0)_{\alpha \beta }=\int _{[\M (\P ^{1},( \alpha, \beta))
]^{vir}}1 =
\begin{cases}
\frac{1}{\combinatfactor (\alpha )}&\text{ if }\alpha =\beta, \\
0& \text{ if } \alpha \neq \beta .
\end{cases}
\]

We use the gluing law to derive the vanishing of
all the terms of $$Z_{d} (0)_{\alpha \beta
}=\tfrac{1}{\combinatfactor (\beta )} Z_{d} (0)^{\beta }_{\alpha }$$ of
higher degree in $t$:
\[
Z_{d} (0)^{\beta }_{\alpha }=Z_{d} (0)^{\gamma }_{\alpha }Z_{d} (0)^{\beta
}_{\gamma }.
\]
For $b>0$ we have
\[
Z_{d}^{b} (0)^{\beta }_{\alpha }=\sum _{\gamma ,\;
b_{1}+b_{2}=b}Z_{d}^{b_{1}} (0)^{\gamma }_{\alpha }Z^{b_{2}}_{d} (0)^{\beta
}_{\gamma }.
\]
By induction, we may assume that $Z_{d}^{b'} (0)^{\gamma }_{\alpha }=0$ for
$0<b'<b$. Then 
\begin{align*}
Z_{d}^{b} (0)^{\beta }_{\alpha }&=\sum _{\gamma }\left(Z_{d}^{b} (0)^{\gamma
}_{\alpha }Z_{d}^{0} (0)^{\beta }_{\gamma }+Z_{d}^{0} (0)^{\gamma
}_{\alpha }Z_{d}^{b} (0)^{\beta }_{\gamma } \right)\\
&=2Z_{d}^{b} (0)^{\beta }_{\alpha }
\end{align*}
and so $Z_{d}^{b} (0)_{\alpha \beta }=0$ for $b>0$. This completes the
proof of Lemma~\ref{lem: local invs of the annulus}. \qed 

All the other axioms of a (1+1)-dimension TQFT given in Section 
\ref{sec: TQFTs} follow immediately from the definitions.
The proof of Proposition~\ref{prop: gluing formulas imply Z
is a TQFT} is complete. 
\qed

To complete the proof of Theorem~\ref{thm: the TQFT exists}, we must
prove $Z_{d} (-)$ is semisimple. By Proposition 
\ref{prop: A/m semi-simple over R/m => A s-s over R},
it suffices to analyze $Z_{d} (-)$ at $t=0$.
Let $$Z^{0}_{d} (-):2\mathbf{Cob}\to \qnums \mathbf{mod}$$ 
denote composition
of $Z_{d} (-)$ with the natural functor $R\mathbf{mod}\to \qnums
\mathbf{mod}$ obtained by setting $t=0$.

\begin{lemma}\label{lemma: Z at t=0 is the DW/FQ TQFT}
$Z^{0}_{d} (-)$ is the
semisimple TQFT (over $\qnums $) given by ``finite gauge theory with gauge
group $S_{d}$'' studied by Dijkgraaf-Witten and Freed-Quinn
\cite{Dijkgraaf-Witten90,Freed-Quinn}. The corresponding Frobenius algebra
is isomorphic to the center of the group algebra $\qnums [S_{d}]$.
\end{lemma}
\textsc{Proof:} 
The invariant $Z^{0}_{d} (g)_{\bfalpha}$ 
is, by definition, the degree of the virtual class of
the expected dimension 0 components of $\M (X,\bfalpha)$. 
Using the dimension formula,  a relative
stable map $$[f:C\to (X,b_{1},\dots ,b_{r})]$$ is easily seen to
lie in an expected dimension
0 component if and only if $f$ is unramified away from the prescribed
ramification points $b_{1},\dots ,b_{r}$. Since such a map has no deformations,
the expected dimension 0 components of $\M (X,\bfalpha)$ 
have actual dimension 0.
Hence, the invariant $Z^{0}_{d} (g)_{\bfalpha}$ 
is equal to a weighted count of maps
with only the prescribed ramification. Each map is
weighted by the reciprocal of the number of automorphisms.
The latter count is, by definition, the Hurwitz number.

Maps with only the prescribed ramification over $b_1, \dots, b_r$
have a gauge theoretic interpretation in terms of principal
$S_{d}$ bundles over an $r$-punctured genus $g$ surface. The 
Hurwitz numbers can be
viewed as counting principal $S_{d}$-bundles over punctured 
surfaces. The associated 
TQFT has been studied in detail and is well-known to
be semisimple \cite{Dijkgraaf-Witten90},
\cite{Freed-Quinn} (but see \cite{Dijkgraaf-mirror95} for a short
explanation). 
The proofs of both Lemma
 \ref{lemma: Z at t=0 is the DW/FQ TQFT} and Theorem 
\ref{thm: the TQFT exists}
are complete.
\qed

The Frobenius algebra $H=\oplus _{\alpha
}\qnums e_{\alpha }$ obtained from $Z_d^0(-)$
is isomorphic to the center of $\qnums [S_{d}]$, the
group ring of the symmetric group. The center has a Frobenius structure
with counit 
\[
\mu \left(\sum _{g\in S_{d}}a_{g}g \right) = \frac{a_{\Id}}{d!}.
\]
The isomorphism is given by
\[
e_{\alpha }\mapsto \sum _{g\in c (\alpha )}g.
\]
The basis $\{e_{\alpha } \}$ is not idempotent, but there is a natural
idempotent basis $\{v_{R} \}$ labeled by irreducible representations $R$ of
$S_{d}$:
\[
v_{R}=\dim R \sum _{\alpha \vdash d}\frac{\chi _{R} (c (\alpha ))}{\combinatfactor
(\alpha )}e_{\alpha }
\]
where $\chi _{R}$ is the character of $R$. 
The calculation,
\[
\mu (v_{R})=\left(\frac{\dim R}{d!} \right)^{2},
\]
yields the following elegant formula counting unramified covers:
\[
Z_{d}^{0} (g)=\sum _{R}\left(\frac{d!}{\dim R} \right)^{2g-2}.
\]

Theorem~\ref{thm: Zd(g) formula} is an immediate consequence of 
Theorem~\ref{thm: the TQFT exists} and the formula
for closed surfaces  (\ref{eqn: formula for
invariant in a semi-simple TQFT}).

\begin{remark}\label{rem: lambda has sqrt - open theory speculation}
We have proven the constant terms of the power series $$\lambda _{1},\dots
,\lambda _{p (d)}$$ appearing in our structure formula (Theorem~\ref{thm:
Zd(g) formula}) are exactly the numbers $$\left(\frac{d!}{\dim
R}\right)^{2}.$$ Hence, the series $\lambda _{i}\in \qnums [[t]]$ have
square roots in $\qnums [[t]]$. The existence of these square roots may
have deeper significance. TQFTs as defined in Section~\ref{sec: TQFTs} are
inherently ``closed string'' TQFTs. One can also axiomatize ``open string''
TQFTs (over a ring $R$) which include the closed theory as a
subsector. Given a closed string TQFT, we may ask what are the possible
open string TQFTs which contain it?  If the closed string TQFT is
semisimple and the corresponding $\lambda_i $'s \emph{have square roots} in
$R$, then there is an elegant classification of the possible open string
TQFT containing the given closed string TQFT. The open string TQFTs are
completely determined by assigning a free $R$-module to each idempotent
basis vector $e_{i}$ and a choice of a sign for the square root of $\lambda
_{i}^{-1}=\mu (e_{i})$. See the lecture notes of G. Moore for a good
discussion \cite{Moore-ITPlecture}.
\end{remark}

\section{Analysis of the TQFT}\label{sec: analysis of TQFT}

In Section \ref{sec: the TQFT}, we constructed $Z_{d} (-)$, a
1+1-dimensional TQFT over $R=\qnums [[t]]$ from the local Gromov-Witten
invariants.  We now analyze the TQFT and the corresponding Frobenius
algebra.

As before, let $H=\oplus _{\alpha \vdash d}Re_{\alpha } $ be
$Z_{d} (S^{1})$. The $R$-module $H$ has the structure of a Frobenius
algebra with multiplication $\cdot $ given by
\[
e_{\alpha }\cdot e_{\beta }=Z_{d} (0)^{\gamma }_{\alpha \beta }e_{\gamma },
\]
unit $\one$ given by 
\[
\one =Z_{d} (0)^{\alpha }e_{\alpha },
\]
comultiplication $\Delta $ given by 
\[
\Delta (e_{\alpha })=Z_{d} (0)^{\beta
\gamma }_{\alpha }e_{\beta }\otimes e_{\gamma },
\]
and counit $\mu $ given by 
\[
\mu (e_{\alpha })=Z_{d} (0)_{\alpha }.
\]

The relative local invariants $Z_{d} (0)_{\alpha }$ and $Z_d(0)_{\alpha
\beta \gamma } $ therefore determine the whole TQFT and hence all the
invariants (note that the invariant $Z_{d} (0)_{\alpha \beta }$ is given by
Lemma~\ref{lem: local invs of the annulus}). The invariants $Z_{d}
(0)_{\alpha \beta \gamma }$, which correspond to the ``pair of pants'' are
in general, difficult to compute.

On the other hand, the invariants $Z_{d} (0)_{\alpha }$ can be computed by
localization.

\begin{theorem}\label{thm: invariants of the disk}
The local relative invariant of $\P ^{1}$ relative to one point is 
\begin{equation}
\label{pqpq}
Z_{d} (0)_{\alpha }= (-1)^{d-l}\frac{t^{d}}{\combinatfactor (\alpha )}\prod
_{i=1}^{l} \left(2\sin (\frac{\alpha _{i}t}{2}) \right)^{-1}
\end{equation}
where the parts of $\alpha \vdash d$ are $\alpha _{1}+\dots +\alpha _{l}=d$.
\end{theorem}

\begin{warning} {\em The local invariants of $\P^1$ relative to one point
were previously studied (in the connected case) by J. Li and Y. Song
\cite{Li-Song}. However, the calculation of \cite{Li-Song} was incomplete
as most localization terms were left unanalyzed by the authors (with the
hope that the contributions vanished).  In fact, the omitted terms of
\cite{Li-Song} {\em do not} vanish, and the calculation there is
wrong. Remarkably, the correct calculation differs only by the sign
$(-1)^{d-l}$. }
\end{warning}

We will calculate the local invariants in the connected case
\begin{equation}
\label{zaza}
\int_{[\overline{M}_g(\P^{1}, \alpha  )]^{vir}} c_b(I(\P^{1}, \alpha)),
\end{equation}
where the relative point is taken to be $\infty$.
The disconnected formula (\ref{pqpq})
will be obtained afterwards by exponentiation.

\begin{lemma} \label{lvann}
The connected local invariants (\ref{zaza}) vanish if $\ell(\alpha)>1$.
\end{lemma}

\textsc{Proof}: We analyze the $K$-theoretic element $I(\P^{1}, \alpha)$,
$$R^{1 }\pi _{*} (f^{*} (\omega_{\mathcal X} (B))\oplus \O (-R))-R^{0 }\pi
_{*} (f^{*} (\omega_{\mathcal X} (B))\oplus \O (-R)).$$ There is a
canonical map, $\epsilon: {\mathcal X} \rightarrow \P^{1},$ obtained by
contracting the destabilizations of the target.  The basic isomorphism,
$$ \omega_{\mathcal X}(B) \stackrel{\sim}{=}
\epsilon^*(\omega_{\P^{1}}(\infty)),$$ is easily checked on the Artin stack of
destabilizations of $\P^{1}$ at $\infty$. Of course, $$\omega_{\P^{1}}(\infty)
\stackrel{\sim}{=} \O_{\P^{1}}(-1).$$ Since both $\O_{\P^{1}}(-1)$ and $\O(-R)$ are
negative, $I({\P^{1}},\alpha)$ simplifies to
\begin{equation}\label{bunn}
R^{1 }\pi _{*} ((\epsilon f)^{*} (\O_{{\P^{1}}} (-1))\oplus
\O (-R)).
\end{equation}
By Riemann-Roch, the rank of the bundle (\ref{bunn}) equals the virtual
dimension $b$ of the moduli space $\overline{M}_{g}({\P^{1}}, \alpha)$.

The sheaf $R^0\pi_*(\mathcal{O}_R)$ is a rank $\ell(\alpha)$ trivial
bundle\footnote{To be precise, $R^{0}\pi_{*} (O_{R})$ is trivial when
pulled back to the finite \'etale cover of $\overline{M} _{g} ({\P^{1}} ,\alpha )
$ given by marking the ramification divisor. It is therefore trivial in
$K$-theory over $\qnums $ and so the following argument holds.} on the
moduli space $\overline{M}_{g}({\P^{1}}, \alpha)$.  The exact sequence
\begin{equation} \label{xxxx}
0 \rightarrow \O^{\ell(\alpha)-1} \rightarrow R^1\pi_*( \O(-R)) 
\rightarrow R^1\pi_*(\O) \rightarrow 0,
\end{equation}
is easily obtained from the ideal sequence
$$0 \rightarrow \O(-R) \rightarrow \O \rightarrow \O_R \rightarrow 0.$$

If $\ell(\alpha)>1$, then the bundle $R^1\pi_*(\O(-R))$ contains
a trivial subfactor. Hence, the local invariant (\ref{zaza}) vanishes
unless $\alpha$ consists of only a single part.
\qed
\vspace{+10pt}

\textsc{Proof of Theorem \ref{thm: invariants of the disk}:}
Let $\alpha$ consist of the single part $a$.
Writing $\O$ as $(\epsilon f)^*(\O_{\P^{1}})$ and using (\ref{xxxx}),
we obtain
\begin{equation}
\label{ksd}
I({\P^{1}},(a)) = R^{1 }\pi _{*} ((\epsilon f)^{*} (\O_{{\P^{1}}} (-1)  \oplus \O_{\P^{1}})) .
\end{equation}
We will use the above form of $I({\P^{1}},(a))$ in the
localization analysis below.

We now define the appropriate
torus actions.
Let ${\P^{1}}={\mathbb P}(V)$ where $$V=\com \oplus \com.$$
Let $\com^*$ act diagonally on $V$:
\begin{equation}
\label{repp}
\xi\cdot (v_1,v_2) = ( v_1, 
\xi \cdot v_2).
\end{equation}
Let $0, \infty$ be the fixed points $(1:0), (0:1)$ of the corresponding
action on ${\P^{1}}$.

The $\com^*$-action on ${\P^{1}}$ canonically lifts to a $\com^*$-action on the
moduli space of maps $\overline{M}_{g}({\P^{1}},(a))$ relative to
$\infty$. Discussions of the virtual localization formula in the relative
context can be found in \cite{Li-Song}, \cite{Faber-Pandharipande-03},
\cite{Graber-Vakil}.

An equivariant lifting  of $\com^*$ to a line bundle $L$
over 
${\P^{1}}$ is uniquely determined by the weights $[l_0,l_\infty]$
of the fiber
representations at the fixed points 
$L_0, L_\infty$.
The canonical lifting of $\com^*$ to the
tangent bundle $T_{\P^{1}}$ has weights $[1,-1]$.
We will utilize the equivariant liftings of
$\com^*$ to $\O_{{\P^{1}}}(-1)$ and $\O_{\P^{1}}$ with weights
$[-1,0]$ and $[0,0]$ respectively.
An equivariant lift of $I({\P^{1}},(a))$
is canonically induced by (\ref{ksd}).

The integral (\ref{zaza}) may now be calculated via the
virtual localization formula. The answer is expressed
as a sum over all localization graphs, see 
\cite{Faber-Pandharipande-03}.
Fortunately,
our choice of equivariant liftings leads to a complete collapse of the sum.

Localization graphs are in bijective correspondence
to the $\com^*$-fixed loci of $\overline{M}_g({\P^{1}},(a))$.
Let $[f]$ be a $\com^*$-fixed point. The vertices of the associated 
graph $\Gamma$ correspond to the connected components of the
set
\begin{equation}
\label{yoyo}
(\epsilon f) ^{-1}( \{0,\infty\}).
\end{equation}
The edges of $\Gamma$ correspond to the components of the domain
of $f$ which map dominantly to ${\P^{1}}$ under $\epsilon f$.

Let $\Gamma$ be a graph with a nonvanishing contribution to the integral
\begin{equation}
\label{zazaza}
\int_{[\overline{M}_g({\P^{1}}, (a)  )]^{vir}} c_b(I({\P^{1}}, (a))),
\end{equation}
Then,
\begin{enumerate}
\item[(i)] Since the monodromy condition $(a)$ over $\infty$ is
            transitive, there can be only 1 connected component of (\ref{yoyo}) over
            $\infty$. Hence, there is a unique vertex $v$ of $\Gamma$ over
     $\infty$. 

\item[(ii)]
The weight $0$ linearization of  $\O_{\P^{1}}(-1)$ over $\infty$ implies 
$v$ has valence 1, see \cite{Fa-Pa}.
\item[(iii)] The vertex $v$ carries the
  class
$c_{g(v)}(\mathbb E^*)^2$ obtained from the
weight 0 linearizations of $\O_{\P^{1}}(-1)$ and $\O_{\P ^{1}}$ over $\infty$. 
As the class  vanishes for $g(v)>0$ by Mumford's relation,
$$c({\mathbb E})\cdot c({\mathbb E^*}) = 1,$$
$v$ must be of genus 0.

\item[(iv)] By (i-iii), $\Gamma$ has a single vertex $v$ of genus 0
 over $\infty$.
The fixed moduli spaces over $\infty$ consist of relative
maps to an {\em unparameterized} bubble with two relative points:
the attaching point of the bubble and the bubbled $\infty$.
 The relative condition over the attaching
point is  $(a)$ by the valence 1 restriction.
The relative condition over the bubbled $\infty$ is $(a)$
by assumption.
The resulting unparameterized moduli space is degenerate, see
\cite{Faber-Pandharipande-03}.
The graph $\Gamma$ must therefore have {\em no contracted components}
over the point $\infty$ of the original ${\P^{1}}$. 
\end{enumerate}
Therefore, a {\em unique} graph $\Gamma$,  consisting
of a genus $g$ vertex over $0$, a unique edge of degree $a$,
and a degenerate genus 0 vertex over $\infty$,
contributes to (\ref{zazaza}).

A straightforward calculation using the virtual localization
formula yields the formula:
\begin{equation}\label{cont}
\int_{[\overline{M}_g({\P^{1}},(a))]^{vir}} c_b(I({\P^{1}},(a)))=
(-1)^{a-1} a^{2g-2}\int_{\overline{M}_{g,1}} \psi^{2g-2}_1 c_g({\mathbb E}).
\end{equation}
We use the terminology of 
\cite{Faber-Pandharipande-03}
for the unique contributing graph $\Gamma$:
\begin{enumerate}
\item[$\bullet$] The automorphism factor $|{\mathbf A}_\Gamma|$ is $a$ obtained from
  the Galois group of the dominant component.
\item[$\bullet$]
The vertex contribution  is
$$\int_{\overline{M}_{g,1}} 
\frac{c(\mathbb E^*) \ (-1)^g c(\mathbb E)\
 (-1)^g  c_g(\mathbb E)}{\frac{1}{a}-\psi},$$
where the last two factors in the numerator are obtained from the
integrand.
\item[$\bullet$] The edge contribution is
$$\frac{(-1)^{a-1} \frac{(a-1)!}{a^{a-1}}} {\frac{a!}{a^a}},$$
where the numerator is obtained from the integrand.
\end{enumerate}
The contribution (\ref{cont}) of $\Gamma$ is given by the product
$$\frac{1}{|{\mathbf A}_\Gamma|} \ \int_{\overline{M}_{g,1}} 
\frac{c(\mathbb E^*)(-1)^g c(\mathbb E) (-1)^g  c_g(\mathbb E)}{\frac{1}{a}-\psi}
\ \ \frac{(-1)^{a-1} \frac{(a-1)!}{a^{a-1}}} {\frac{a!}{a^a}}$$
after simplification using Mumford's relation.

Using the Hodge integral computation of \cite{Fa-Pa}, we find
\begin{equation*}
\sum_{g\geq 0} t^{b}
\int_{[\overline{M}_g({\P^{1}},(a))]^{vir}} c_b(I({\P^{1}},(a)))=
{(-1)^{a-1}}\frac{t^{a}}{a }
\left( 2\sin(\frac{a t}{2})\right)^{-1}.
\end{equation*}
The proof of Theorem \ref{thm: invariants of the disk} is completed by
taking the associated disconnected integrals.  \qed

\vspace{10pt}

As a corollary to Theorem~\ref{thm: invariants of the disk} we can
completely determine the TQFT $Z_{1} (-)$ and hence all the $d=1$ relative
invariants.
\begin{cor}\label{cor: d=1 full TQFT}
The TQFT $Z_{1} (-)$ is isomorphic to $R_{\lambda }$ where 
\[
\lambda =\left(\frac{\sin (t/2)}{t/2} \right)^{2}
\]
via the isomorphism
\[
e_{[1]}\mapsto \left( \frac{\sin (t/2)}{t/2} \right)^{-1}.
\]
Consequently, the $d=1$ relative invariants are given by
\[
Z_{1} (g)_{\alpha _{1},\dots ,\alpha _{r}}=\left(\frac{\sin (t/2)}{t/2}
\right)^{2g-2+r}.
\]
\end{cor}

\appendix \section{The proof of the gluing formulas}\label{appendix: proof of
gluing}

To prove the gluing formulas in Theorem~\ref{thm: gluing laws}, we will
consider
algebraic degenerations corresponding to the (topological) splittings of
the TQFT. To simplify the exposition, we first derive Equation~(\ref{eqn:
gluing along separating curve}) of Theorem~\ref{thm: gluing laws} with $r=0$.

Consider the nodal curve 
\[
W_{0}=X_{1}\bigcup _{b_{1}=b_{2}}X_{2}
\]
obtained by joining non-singular curves $X_{1} $ and $X_{2}$ of genus
$g_{1} $ and $g_{2}$ at points $b_{i}\in X_{i}$. Let $W\to \affineline $ be
a generic, 1-parameter deformation of $W_{0}$ for which the fibers $W_{t}$
for $t\neq 0\in \affineline $ are nonsingular curves of genus
$g=g_{1}+g_{2}$. (The base of the degeneration can be any smooth curve. For
simplicity, we take it to be $\affineline $.)

In Li's theory, the moduli of relative stable maps spaces arise by constructing a
good limit for the moduli spaces $\M (W_{t} )$ as $t$ approaches 0. Li's method
involves
a stack $\W$ of expanded degenerations of $W$: an
Artin stack over $\affineline $ which, in addition to $W$, includes
degenerations to the curves $W[n]_{0}$ obtained by inserting a chains of
$\P ^{1}$'s between $X_{1}$ and $X_{2}$:
\[
W[n]_{0}=X_{1}\cup \underbrace{\P ^{1}\cup \dots \cup \P ^{1}}_{n-1}\cup X_{2}.
\]

Following Li (but with our conventions (i-iii) of Definition~\ref{defn: rel
stable maps} regarding the domain curve), we define $\M (\W)$ to be the
stack of non-degenerate, pre-deformable, degree $d$ stable maps to $\W$
(see \cite{Li-relative1} section~3 for the definitions of non-degenerate,
pre-deformable, and $\W$). Li proves that $\M (\W )$ is a Deligne-Mumford
stack (Theorem~0.1 of \cite{Li-relative1}). Each component of $\M(\W)$ of
given fixed expected dimension is proper and separated over $\affineline $
(also Theorem~0.1 of \cite{Li-relative1}).

$\M (\W_{0} )$, the central fiber of $\M (\W )$, can be expressed, up to
finite covers, as the union of products of relative stable map moduli
spaces. Moreover, $\M (\W )$ has a virtual fundamental class whose
intersection with $\M (\W_{t} )$ for $t\neq 0$ is the usual virtual
fundamental class $[\M (W_{t} )]^{vir}$, and whose intersection with $\M
(\W _{0} )$ is compatible with the decomposition into relative stable map
spaces. To be precise, Li's virtual cycle formula (Theorem~3.15 of
\cite{Li-relative2}), adapted to our setting and conventions is:
\begin{equation}\label{eqn: virtual cycle formula}
[\M (\W _{0} )]^{vir}=\sum _{\alpha \vdash d}\combinatfactor (\alpha ) \left(\Phi
_{\alpha } \right)_{*} \left([\M (X_{1},\alpha) ]^{vir}\times 
[\M (X_{2},{\alpha })]^{vir} \right).
\end{equation}
Here the map
\[
\Phi _{\alpha }:\M (X_{1},{\alpha })\times \M (X_{2}, {\alpha })\to \M (\W
_{0} )
\]
is obtained by constructing a family of $\M (\W _{0} )$ maps from the
universal maps over $\M (X_{i},{\alpha })$ by gluing along the universal
prescribed ramification and branch divisors:

\begin{equation}\label{diagram: universal diagram for Phi}
\begin{diagram}[height=1.5cm]
U_{1}\bigcup
_{R_{1}=R_{2}}U_{2}&\rTo^{f_{1}\cup
f_{2}\quad \quad \quad }&\X_{1}\bigcup
_{B_{1}=B_{2}}\X_{2}\\
\dTo^{\pi _{1}\cup \pi _{2}}&\ldTo_{p_{1}\cup p_{2}}&\\
\M (X_{1},{\alpha })\times \M (X_{2},{\alpha})
&&
\end{diagram}
\end{equation}

Strictly speaking, we must to pass to the finite \'etale cover of $$\M
(X_{1},{\alpha })\times \M (X_{2},{\alpha })$$ obtained by marking the
ramification divisor. The ordering is necessary to obtain the
identification $R_{1}=R_{2}$.  However, the map $(\Phi _{\alpha })_{*}$ on
$\qnums $-cycles is well defined.  The degree of the finite \'etale map is
included in our constant in the virtual cycle formula. Since our target is
a curve and the branch divisor is a point, the diagonal constraint which
occurs in Li's general cycle formula does not appear here.

In order to apply the virtual cycle formula to obtain our gluing formulas,
we will define a $K$-theory class $I (\W )$ on $\M (\W )$ with the
following properties:
\begin{enumerate}
\item[(i)]
For $t\neq
0$, $I(\W)$ restricts to $I (W_{t})$ on $\M (W_{t} )$.
\item[(ii)]
For $t=0$, $I(\W)$ restricts to
a class that pulls back via $\Phi _{\alpha }$ to $$I (X_{1},{\alpha })\oplus
I (X_{2},{\alpha })$$ on $\M (X_{1},{\alpha })\times \M (X_{2},{\alpha})$. 
\end{enumerate}

\begin{proposition}\label{prop: I(W) at 0 pulls back to the rel I}
Let $\pi :U\to \M (\W )$, $p:\X\to \M (\W )$, and
$f:U\to \X$ be the universal domain, universal target,
and the universal map for $\M (\W )$. Let $\omega _{p}$ be the
relative dualizing sheaf of the universal target. Define the $K$-theory
class $I (\W )$ by
\[
I (\W )=-R^{\bullet }\pi _{*}f^{*} (\omega _{p}\oplus
\O_{\X}).
\]
Let $I (\W )_{t}$ be the restriction of $I (\W )$ to $\M (\W _{t})$. Then
for $t\neq 0$,
\[
I (\W )_{t}=I (W_{t})
\]
and for $t=0$ we have
\[
\Phi _{\alpha }^{*} (I (\W )_{0})=I (X_{1},{\alpha })\oplus I
(X_{2},{\alpha }).
\]
\end{proposition}

\textsc{Proof:} First, the universal family $p:\X\to \M (\W )$ is a flat
family of prestable curves, so there exists a relative dualizing sheaf 
$\ \omega
_{p}$. Over $t\neq 0$, the family is constant with fiber $W_{t}$,
so $$\omega _{p_t} \stackrel{\sim}{=}K_{W_{t}},$$ and thus $I (\W
)_{t}=I (W_{t})$.

We now compute $\Phi _{\alpha }^{*} (I (\W )_{0})$. Let $\pi _{\alpha
}=\pi _{1}\cup \pi _{2}$, $f_{\alpha }=f_{1}\cup f_{2}$, and $p_{\alpha
}=p_{1}\cup p_{2}$ denote the maps in diagram~(\ref{diagram: universal
diagram for Phi}), and let $U_{\alpha }=U_{1}\cup
U_{2}$ and $\X_{\alpha }=\X_{1}\cup
\X_{2}$. By the definition of $\Phi _{\alpha }$, we have
\[
\Phi _{\alpha }^{*} (I (\W )_{0})=-R^{\bullet } (\pi _{\alpha })_{*}
f_{\alpha }^{*} (\omega _{p_{\alpha }}\oplus \O_{\X_{\alpha}}).
\]

Consider the following two short exact sequences of sheaves on $U_{\alpha }$.
\begin{align*}
0\to \O_{U_{\alpha }}\longrightarrow \O_{U_{1}}&\oplus
\O_{U_{2}}\longrightarrow \O_{R}\to 0\\
0\to f_\alpha^*(\omega _{p_{\alpha }})\to f_1^*(\omega _{p_{1}} (B_{1}))
&\oplus
f_2^*(\omega _{p_{2}} (B_{2}))\to \O_{R}\to 0.
\end{align*}
The first sequence is the usual normalization sequence. The second is
obtained from
standard facts about the dualizing sheaf $\omega $ of a nodal
curve. 
We obtain the
following equalities in $K$-theory:
$$
f_{\alpha }^{*} (\O _{\X_{\alpha }})=\O_{U_{1}}+\O_{U_{2}}-\O_{R}, $$
$$
f_{\alpha }^{*} (\omega _{p_{\alpha }})=f_{1}^{*} (\omega _{p_{1}}
(B_{1}))+f_{2}^{*} (\omega _{p_{2}}
(B_{2})) - \O_{R},$$
where the first equation uses the isomorphism 
$\O_{U_\alpha} \stackrel{\sim}{=}f_{\alpha }^{*} (\O _{\X_{\alpha }})$.

The divisor sequence for $R_{i}\subset U_{i}$ yields
\[
\O_{U_{i}}=\O_{R}+ \O_{U_{i}} (-R_{i})
\]
in $K$-theory, so
\[
f_{\alpha }^{*} (\omega _{p_{\alpha }}\oplus \mathcal{O_{X_{\alpha
}}})=\sum _{i=1}^{2}f^{*}_{i} (\omega _{p_{i}}
(B_{i}))+\sum _{i=1}^{2}\O_{U_{i}} (-R_{i})
.
\]
After
applying $-R^{\bullet } (\pi _{\alpha })_{*}$ to both sides, we obtain
\begin{align*}
\Phi _{\alpha }^{*} (I (\W )_{0})&=\sum _{i=1}^{2}-R^{\bullet } (\pi
_{i})_{*}\left(f^{*}_{i} (\omega _{p_{i}}
(B_{i}))+\O_{U_{i}} (-R_{i})\right)\\
&=\sum _{i=1}^{2}I (X_{i},{\alpha }).
\end{align*}
The proof of the Proposition is complete. \qed
\vspace{10pt}

With the class $I (\W )$ in hand, the gluing formula is proven as
follows. Using the fact that $\int _{[\M (\W _{t})]^{vir}}c_{b} (I (\W ))$
is independent of $t$ we compute:

\begin{eqnarray*}
Z^{b}_{d} (g_{1}+g_{2})&=&\int   _{[\M (\W_{t} )]^{vir}}c_{b} 
(I(\W )_{t} )\\
&&\\
&=&\int _{[\M (\W _{0})]^{vir}}c_{b} (I (\W )_{0})\\
&\quad &\\
&=&\quad \sum _{\alpha \vdash d}\quad \combinatfactor (\alpha )\; \int   _{[\M (X_{1}
,{\alpha })]^{vir}\times [\M (X_{2},{\alpha })]^{vir} }c_{b}\left(I (X_{1},{\alpha })\oplus I (X_{2},{\alpha }) \right)\\
&&\\
&=&\sum _{ \begin{smallmatrix}{\alpha \vdash d}\\
 {b_{1}+b_{2}=b}  \end{smallmatrix}}
\combinatfactor (\alpha ) \int
  _{[\M (X_{1},{\alpha })]^{vir}}c_{b_{1}} (I (X_{1},{\alpha }))
\int _{[\M (X_{2},{\alpha })]^{vir}}c_{b_{2}} (I
(X_{2},{\alpha }))\\
&=&\sum _{ \begin{smallmatrix}{\alpha \vdash d}\\
 {b_{1}+b_{2}=b}  \end{smallmatrix}}
 \combinatfactor (\alpha )\;
Z_{d}^{b_{1}} (g_{1})_{\alpha }\; Z^{b_{2}}_{d} (g_{2})_{\alpha }.
\end{eqnarray*}
We have proven the gluing formula (\ref{eqn: gluing along separating
curve}) in case  $r=0$. 

The proof of the second gluing
formula (\ref{eqn: gluing along non-separating curve}) for $r=0$
 is almost identical. We consider a degeneration
\[
W\to \affineline 
\]
where $W_{t}$ for $t\neq 0$ is a nonsingular genus $(g+1)$ curve and 
\[
W_{0}=X/b_{1}\sim b_{2}
\]
is an irreducible nodal curve whose normalization $(X,b_{1},b_{2})$ is a
smooth genus $g$ curve with marked points $b_{1},b_{2}\in X$ lying over the
node.

As in the previous case, we construct a stack $\W $ of expanded
degenerations of $W$ and define $\M (\W )$ to be the stack of non-degenerate,
pre-deformable, stable maps to $\W $. The cycle formula is now 
\[
[\M (\W _{0})]^{vir}=\sum _{\alpha \vdash d}\combinatfactor (\alpha ) (\Phi
_{\alpha \alpha })_{*} ([\M (X, (\alpha, \alpha))]^{vir})
\]
where 
\begin{equation}
\Phi _{\alpha \alpha } :\M (X,({\alpha, \alpha }))\to \M (\W _{0})
\end{equation}
is obtained by gluing together the two universal prescribed branched
divisors and the two universal prescribed ramification divisors over $ \M
(X,({\alpha, \alpha }))$:
\[
\begin{diagram}[height=1cm]
R_{1},R_{2} \rInto &U_{\alpha \alpha }&\rTo^{f_{\alpha \alpha }}&\X _{\alpha \alpha } &\lInto B_{1},B_{2}\\
&
\dTo^{n}&&\dTo_{n}&\\
\quad \quad R \rInto &U&\rTo^{f}&\X &\lInto B \quad \quad \\
&\dTo^{\pi }&\ldTo_{p}&&\\
&\M (X,({\alpha, \alpha }))&
\end{diagram}
\]

In the above diagram, $\pi _{\alpha \alpha }=\pi \circ n$ and $p_{\alpha \alpha
}=p\circ n$ are the universal domain and universal range for $\M
(X,({\alpha, \alpha }))$ respectively and the stacks
\[
\X =\X _{\alpha \alpha }/B_{1}\sim B_{2} \quad \quad U=U_{\alpha \alpha
}/R_{1}\sim R_{2}
\]
are obtained by gluing together the two universal prescribed branched
divisors and the two universal prescribed ramification divisors
respectively (after possibly passing to an \'etale cover). The $\M (\W
_{0})$ family given by $(\pi ,f,p)$ defines the map $\Phi _{\alpha \alpha
}$.

As before, we define $I (\W )=-R^{\bullet }\pi _{*}f^{*} (\omega _{p}\oplus
\O_{\X } )$. The analogue of Proposition~\ref{prop: I(W) at 0 pulls back to
the rel I} is the assertion
\begin{equation}
\label{asss}
\Phi ^{*}_{\alpha \alpha } (I (\W )_{0})=I (X,({\alpha ,\alpha })). 
\end{equation}
which is proven in essentially the same way:

By applying $R^{\bullet }\pi _{*} (-)$ to the two exact sequences
\[
0\to \O _{U} \to n_{*}\O _{U_{\alpha \alpha }}\to \O _{R}\to 0
\]

\[
0\to f^{*}\omega _{p}\to n_{*}f^{*}_{\alpha \alpha } (\omega _{p_{\alpha \alpha }} (B_{1}+B_{2}))\to \O _{R}\to 0
\]
we obtain the following equalities in the $K$-theory of $\M (X,(\alpha
,\alpha ) )$:
\begin{align*}
R^{\bullet }\pi _{*}\O _{U}&=R^{\bullet }\pi _{\alpha \alpha *}\O _{U_{\alpha \alpha }} - R^{\bullet }\pi _{*}\O _{R}\\
R^{\bullet }\pi _{*}f^{*}\omega _{p}&=R^{\bullet }\pi _{\alpha \alpha *}f^{*}_{\alpha \alpha } (\omega _{p_{\alpha \alpha }} (B_{1}+B_{2}))-R^{\bullet }\pi _{*}\O _{R}.
\end{align*}

Therefore we get
\begin{align}\label{eqn: expression for Phiaa(I(W))}
\Phi ^{*}_{\alpha \alpha } (I (\W )_{0})&=-R^{\bullet }\pi _{*}f^{*}
(\omega _{p}\oplus \O _{\X })\nonumber\\
&=-R^{\bullet }\pi _{\alpha \alpha *}f^{*}_{\alpha \alpha } (\omega _{p_{\alpha \alpha }} (B_{1}+B_{2}))+R^{\bullet }\pi _{*}\O _{R}\nonumber \\
&\quad -R^{\bullet }\pi _{\alpha \alpha *}\O _{U_{\alpha \alpha }}\quad
\quad \quad \quad \quad \quad +R^{\bullet }\pi _{*}\O _{R}.
\end{align}

Since $n_{*}\O _{R_{i}}=\O _{R}$ we get the following equalities in
$K$-theory:
\begin{align*}
2R^{\bullet }\pi _{*}\O _{R}-R^{\bullet }\pi
_{\alpha \alpha *}\O _{U_{\alpha \alpha }}&=R^{\bullet }\pi _{\alpha \alpha
*}\O _{R_{1}}+R^{\bullet }\pi _{\alpha \alpha *}\O _{R_{2}}-R^{\bullet }\pi
_{\alpha \alpha *}\O _{U_{\alpha \alpha }}\\
&=-R^{\bullet }\pi _{\alpha \alpha *}\O _{U_{\alpha \alpha }} (-R_{1}-R_{2})
\end{align*}
where the last equality comes from the divisor sequence for $R_{1}+R_{2}$.

Combining the above with equation \ref{eqn: expression for Phiaa(I(W))}, we
get
\begin{align*}
\Phi ^{*}_{\alpha \alpha } (I (\W )_{0})&=-R^{\bullet }\pi _{\alpha \alpha *} (f^{*}_{\alpha \alpha } (\omega _{p_{\alpha \alpha }} (B_{1}+B_{2}))\oplus \O _{U_{\alpha \alpha }} (-R_{1}-R_{2}))\\
&=I (X,(\alpha ,\alpha ))
\end{align*}
which completes the proof of Equation (\ref{asss}).

The derivation of the gluing formula from Equation (\ref{asss}) is
easily obtained:

\begin{align*}
Z^{b}_{d} (g+1)&=\int _{[\M (\W _{t})]^{vir}}c_{b} (I (\W )_{t})\\
&=\int _{[\M (\W _{0})]^{vir}}c_{b} (I (\W )_{0})\\
&=\sum _{\alpha \vdash d}\combinatfactor (\alpha )\int _{[\M (X,({\alpha, \alpha
})]^{vir}}c_{b} (I (X,({\alpha, \alpha }))\\
&=\sum _{\alpha \vdash d}\combinatfactor (\alpha )Z^{b}_{d} (g)_{\alpha
\alpha }.
\end{align*}

The proof of the gluing formulas in case $r>0$ is identical.
Although Li does not explicitly state the virtual cycle formula
necessary for the $r>0$ case, the techniques and results of 
\cite{Li-relative2}
extend in a
straightforward way. \qed

\subsubsection{Acknowledgments} The authors warmly thank Domenico Fiorenza,
Jun Li, Andrei Okounkov, Michael Thaddeus, and Ravi Vakil for helpful
discussions. Jim Bryan is supported by NSERC, NSF, and the Sloan
Foundation; Rahul Pandharipande is supported by NSF and the Sloan and
Packard foundations.


\end{document}

%% file: cob4.eepic
\setlength{\unitlength}{0.00053333in}
\begingroup\makeatletter\ifx\SetFigFont\undefined%
\gdef\SetFigFont#1#2#3#4#5{%
  \reset@font\fontsize{#1}{#2pt}%
  \fontfamily{#3}\fontseries{#4}\fontshape{#5}%
  \selectfont}%
\fi\endgroup%
{\renewcommand{\dashlinestretch}{30}
\begin{picture}(7396,3865)(0,-10)
\Thicklines
\path(525,3040)(526,3040)(530,3037)
	(539,3032)(552,3024)(569,3014)
	(586,3003)(602,2992)(616,2981)
	(627,2971)(636,2961)(644,2951)
	(650,2940)(655,2928)(660,2915)
	(663,2901)(666,2885)(668,2869)
	(668,2852)(668,2836)(666,2820)
	(663,2804)(660,2790)(655,2777)
	(650,2765)(644,2754)(636,2744)
	(627,2734)(616,2724)(602,2713)
	(586,2702)(569,2691)(552,2681)
	(539,2673)(530,2668)(526,2665)(525,2665)
\Thicklines
\path(1425,790)(1426,792)(1429,795)
	(1435,801)(1442,810)(1453,822)
	(1466,836)(1482,852)(1500,870)
	(1519,888)(1539,907)(1561,925)
	(1584,942)(1609,958)(1635,973)
	(1663,987)(1694,998)(1727,1007)
	(1763,1013)(1800,1015)(1837,1013)
	(1873,1007)(1906,998)(1937,987)
	(1965,973)(1991,958)(2016,942)
	(2039,925)(2061,907)(2081,888)
	(2100,870)(2118,852)(2134,836)
	(2147,822)(2158,810)(2165,801)
	(2171,795)(2174,792)(2175,790)
\Thicklines
\path(1341,2377)(1341,2376)(1342,2373)
	(1344,2368)(1346,2360)(1350,2351)
	(1355,2340)(1361,2328)(1368,2316)
	(1377,2303)(1388,2289)(1403,2275)
	(1420,2261)(1442,2246)(1469,2230)
	(1500,2215)(1530,2202)(1560,2192)
	(1586,2183)(1609,2175)(1627,2170)
	(1641,2166)(1652,2162)(1661,2160)
	(1670,2158)(1678,2156)(1687,2154)
	(1700,2152)(1716,2150)(1737,2148)
	(1764,2145)(1796,2143)(1834,2141)
	(1875,2140)(1916,2141)(1954,2143)
	(1986,2145)(2013,2148)(2035,2150)
	(2051,2152)(2063,2154)(2073,2156)
	(2082,2158)(2090,2160)(2100,2162)
	(2111,2166)(2125,2170)(2143,2175)
	(2165,2183)(2192,2192)(2221,2202)
	(2250,2215)(2280,2230)(2305,2246)
	(2326,2261)(2342,2275)(2354,2289)
	(2364,2303)(2372,2316)(2378,2328)
	(2382,2340)(2385,2351)(2388,2360)
	(2389,2368)(2390,2373)(2391,2376)(2391,2377)
\path(1500,2215)(1501,2217)(1504,2220)
	(1510,2226)(1517,2235)(1528,2247)
	(1541,2261)(1557,2277)(1575,2295)
	(1594,2313)(1614,2332)(1636,2350)
	(1659,2367)(1684,2383)(1710,2398)
	(1738,2412)(1769,2423)(1802,2432)
	(1838,2438)(1875,2440)(1912,2438)
	(1948,2432)(1981,2423)(2012,2412)
	(2040,2398)(2066,2383)(2091,2367)
	(2114,2350)(2136,2332)(2156,2313)
	(2175,2295)(2193,2277)(2209,2261)
	(2222,2247)(2233,2235)(2240,2226)
	(2246,2220)(2249,2217)(2250,2215)
\Thicklines
\path(1341,3202)(1341,3201)(1342,3198)
	(1344,3193)(1346,3185)(1350,3176)
	(1355,3165)(1361,3153)(1368,3141)
	(1377,3128)(1388,3114)(1403,3100)
	(1420,3086)(1442,3071)(1469,3055)
	(1500,3040)(1530,3027)(1560,3017)
	(1586,3008)(1609,3000)(1627,2995)
	(1641,2991)(1652,2987)(1661,2985)
	(1670,2983)(1678,2981)(1687,2979)
	(1700,2977)(1716,2975)(1737,2973)
	(1764,2970)(1796,2968)(1834,2966)
	(1875,2965)(1916,2966)(1954,2968)
	(1986,2970)(2013,2973)(2035,2975)
	(2051,2977)(2063,2979)(2073,2981)
	(2082,2983)(2090,2985)(2100,2987)
	(2111,2991)(2125,2995)(2143,3000)
	(2165,3008)(2192,3017)(2221,3027)
	(2250,3040)(2280,3055)(2305,3071)
	(2326,3086)(2342,3100)(2354,3114)
	(2364,3128)(2372,3141)(2378,3153)
	(2382,3165)(2385,3176)(2388,3185)
	(2389,3193)(2390,3198)(2391,3201)(2391,3202)
\path(1500,3040)(1501,3042)(1504,3045)
	(1510,3051)(1517,3060)(1528,3072)
	(1541,3086)(1557,3102)(1575,3120)
	(1594,3138)(1614,3157)(1636,3175)
	(1659,3192)(1684,3208)(1710,3223)
	(1738,3237)(1769,3248)(1802,3257)
	(1838,3263)(1875,3265)(1912,3263)
	(1948,3257)(1981,3248)(2012,3237)
	(2040,3223)(2066,3208)(2091,3192)
	(2114,3175)(2136,3157)(2156,3138)
	(2175,3120)(2193,3102)(2209,3086)
	(2222,3072)(2233,3060)(2240,3051)
	(2246,3045)(2249,3042)(2250,3040)
\path(2766,3352)(2766,3351)(2767,3348)
	(2769,3343)(2771,3335)(2775,3326)
	(2780,3315)(2786,3303)(2793,3291)
	(2802,3278)(2813,3264)(2828,3250)
	(2845,3236)(2867,3221)(2894,3205)
	(2925,3190)(2955,3177)(2985,3167)
	(3011,3158)(3034,3150)(3052,3145)
	(3066,3141)(3077,3137)(3086,3135)
	(3095,3133)(3103,3131)(3112,3129)
	(3125,3127)(3141,3125)(3162,3123)
	(3189,3120)(3221,3118)(3259,3116)
	(3300,3115)(3341,3116)(3379,3118)
	(3411,3120)(3438,3123)(3460,3125)
	(3476,3127)(3488,3129)(3498,3131)
	(3507,3133)(3515,3135)(3525,3137)
	(3536,3141)(3550,3145)(3568,3150)
	(3590,3158)(3617,3167)(3646,3177)
	(3675,3190)(3705,3205)(3730,3221)
	(3751,3236)(3767,3250)(3779,3264)
	(3789,3278)(3797,3291)(3803,3303)
	(3807,3315)(3810,3326)(3813,3335)
	(3814,3343)(3815,3348)(3816,3351)(3816,3352)
\path(2925,3190)(2926,3192)(2929,3195)
	(2935,3201)(2942,3210)(2953,3222)
	(2966,3236)(2982,3252)(3000,3270)
	(3019,3288)(3039,3307)(3061,3325)
	(3084,3342)(3109,3358)(3135,3373)
	(3163,3387)(3194,3398)(3227,3407)
	(3263,3413)(3300,3415)(3337,3413)
	(3373,3407)(3406,3398)(3437,3387)
	(3465,3373)(3491,3358)(3516,3342)
	(3539,3325)(3561,3307)(3581,3288)
	(3600,3270)(3618,3252)(3634,3236)
	(3647,3222)(3658,3210)(3665,3201)
	(3671,3195)(3674,3192)(3675,3190)
\path(1266,952)(1266,951)(1267,948)
	(1269,943)(1271,935)(1275,926)
	(1280,915)(1286,903)(1293,891)
	(1302,878)(1313,864)(1328,850)
	(1345,836)(1367,821)(1394,805)
	(1425,790)(1455,777)(1485,767)
	(1511,758)(1534,750)(1552,745)
	(1566,741)(1577,737)(1586,735)
	(1595,733)(1603,731)(1612,729)
	(1625,727)(1641,725)(1662,723)
	(1689,720)(1721,718)(1759,716)
	(1800,715)(1841,716)(1879,718)
	(1911,720)(1938,723)(1960,725)
	(1976,727)(1988,729)(1998,731)
	(2007,733)(2015,735)(2025,737)
	(2036,741)(2050,745)(2068,750)
	(2090,758)(2117,767)(2146,777)
	(2175,790)(2205,805)(2230,821)
	(2251,836)(2267,850)(2279,864)
	(2289,878)(2297,891)(2303,903)
	(2307,915)(2310,926)(2313,935)
	(2314,943)(2315,948)(2316,951)(2316,952)
\path(4566,2902)(4566,2901)(4567,2898)
	(4569,2893)(4571,2885)(4575,2876)
	(4580,2865)(4586,2853)(4593,2841)
	(4602,2828)(4613,2814)(4628,2800)
	(4645,2786)(4667,2771)(4694,2755)
	(4725,2740)(4755,2727)(4785,2717)
	(4811,2708)(4834,2700)(4852,2695)
	(4866,2691)(4877,2687)(4886,2685)
	(4895,2683)(4903,2681)(4912,2679)
	(4925,2677)(4941,2675)(4962,2673)
	(4989,2670)(5021,2668)(5059,2666)
	(5100,2665)(5141,2666)(5179,2668)
	(5211,2670)(5238,2673)(5260,2675)
	(5276,2677)(5288,2679)(5298,2681)
	(5307,2683)(5315,2685)(5325,2687)
	(5336,2691)(5350,2695)(5368,2700)
	(5390,2708)(5417,2717)(5446,2727)
	(5475,2740)(5505,2755)(5530,2771)
	(5551,2786)(5567,2800)(5579,2814)
	(5589,2828)(5597,2841)(5603,2853)
	(5607,2865)(5610,2876)(5613,2885)
	(5614,2893)(5615,2898)(5616,2901)(5616,2902)
\path(4566,1177)(4566,1176)(4567,1173)
	(4569,1168)(4571,1160)(4575,1151)
	(4580,1140)(4586,1128)(4593,1116)
	(4602,1103)(4613,1089)(4628,1075)
	(4645,1061)(4667,1046)(4694,1030)
	(4725,1015)(4755,1002)(4785,992)
	(4811,983)(4834,975)(4852,970)
	(4866,966)(4877,962)(4886,960)
	(4895,958)(4903,956)(4912,954)
	(4925,952)(4941,950)(4962,948)
	(4989,945)(5021,943)(5059,941)
	(5100,940)(5141,941)(5179,943)
	(5211,945)(5238,948)(5260,950)
	(5276,952)(5288,954)(5298,956)
	(5307,958)(5315,960)(5325,962)
	(5336,966)(5350,970)(5368,975)
	(5390,983)(5417,992)(5446,1002)
	(5475,1015)(5505,1030)(5530,1046)
	(5551,1061)(5567,1075)(5579,1089)
	(5589,1103)(5597,1116)(5603,1128)
	(5607,1140)(5610,1151)(5613,1160)
	(5614,1168)(5615,1173)(5616,1176)(5616,1177)
\path(4725,1015)(4726,1017)(4729,1020)
	(4735,1026)(4742,1035)(4753,1047)
	(4766,1061)(4782,1077)(4800,1095)
	(4819,1113)(4839,1132)(4861,1150)
	(4884,1167)(4909,1183)(4935,1198)
	(4963,1212)(4994,1223)(5027,1232)
	(5063,1238)(5100,1240)(5137,1238)
	(5173,1232)(5206,1223)(5237,1212)
	(5265,1198)(5291,1183)(5316,1167)
	(5339,1150)(5361,1132)(5381,1113)
	(5400,1095)(5418,1077)(5434,1061)
	(5447,1047)(5458,1035)(5465,1026)
	(5471,1020)(5474,1017)(5475,1015)
\path(2766,652)(2766,651)(2767,648)
	(2769,643)(2771,635)(2775,626)
	(2780,615)(2786,603)(2793,591)
	(2802,578)(2813,564)(2828,550)
	(2845,536)(2867,521)(2894,505)
	(2925,490)(2955,477)(2985,467)
	(3011,458)(3034,450)(3052,445)
	(3066,441)(3077,437)(3086,435)
	(3095,433)(3103,431)(3112,429)
	(3125,427)(3141,425)(3162,423)
	(3189,420)(3221,418)(3259,416)
	(3300,415)(3341,416)(3379,418)
	(3411,420)(3438,423)(3460,425)
	(3476,427)(3488,429)(3498,431)
	(3507,433)(3515,435)(3525,437)
	(3536,441)(3550,445)(3568,450)
	(3590,458)(3617,467)(3646,477)
	(3675,490)(3705,505)(3730,521)
	(3751,536)(3767,550)(3779,564)
	(3789,578)(3797,591)(3803,603)
	(3807,615)(3810,626)(3813,635)
	(3814,643)(3815,648)(3816,651)(3816,652)
\path(2925,490)(2926,492)(2929,495)
	(2935,501)(2942,510)(2953,522)
	(2966,536)(2982,552)(3000,570)
	(3019,588)(3039,607)(3061,625)
	(3084,642)(3109,658)(3135,673)
	(3163,687)(3194,698)(3227,707)
	(3263,713)(3300,715)(3337,713)
	(3373,707)(3406,698)(3437,687)
	(3465,673)(3491,658)(3516,642)
	(3539,625)(3561,607)(3581,588)
	(3600,570)(3618,552)(3634,536)
	(3647,522)(3658,510)(3665,501)
	(3671,495)(3674,492)(3675,490)
\thinlines
\dottedline{90}(150,1390)(150,790)
\dottedline{90}(150,3490)(150,2890)
\dottedline{90}(6825,3415)(6825,2815)
\Thicklines
\path(4725,2740)(4726,2742)(4729,2745)
	(4735,2751)(4742,2760)(4753,2772)
	(4766,2786)(4782,2802)(4800,2820)
	(4819,2838)(4839,2857)(4861,2875)
	(4884,2892)(4909,2908)(4935,2923)
	(4963,2937)(4994,2948)(5027,2957)
	(5063,2963)(5100,2965)(5137,2963)
	(5173,2957)(5206,2948)(5237,2937)
	(5265,2923)(5291,2908)(5316,2892)
	(5339,2875)(5361,2857)(5381,2838)
	(5400,2820)(5418,2802)(5434,2786)
	(5447,2772)(5458,2760)(5465,2751)
	(5471,2745)(5474,2742)(5475,2740)
\thinlines
\dottedline{90}(525,2215)(525,1915)
\dottedline{90}(6825,1165)(6825,565)
\dottedline{90}(6450,2140)(6450,1840)
\dottedline{90}(2775,2215)(2775,1615)
\dottedline{90}(3825,2215)(3825,1615)
\dottedline{90}(5100,2290)(5100,1540)
\Thicklines
\path(525,340)(526,339)(528,339)
	(530,337)(534,335)(539,333)
	(546,330)(555,327)(565,323)
	(577,318)(591,313)(607,308)
	(626,302)(647,297)(670,291)
	(697,284)(727,278)(762,271)
	(800,264)(844,256)(894,248)
	(950,239)(1014,231)(1086,221)
	(1166,211)(1254,201)(1350,190)
	(1414,183)(1479,176)(1544,170)
	(1607,163)(1668,157)(1726,151)
	(1781,145)(1833,139)(1880,134)
	(1924,128)(1964,123)(2001,118)
	(2034,114)(2063,109)(2090,105)
	(2114,101)(2136,97)(2156,93)
	(2174,89)(2192,85)(2209,81)
	(2225,77)(2242,74)(2259,70)
	(2277,67)(2297,64)(2319,60)
	(2343,57)(2370,54)(2400,51)
	(2434,48)(2472,45)(2515,43)
	(2562,40)(2614,38)(2672,36)
	(2734,35)(2802,34)(2876,33)
	(2953,33)(3036,34)(3121,35)
	(3210,37)(3300,40)(3390,44)
	(3479,48)(3565,54)(3648,60)
	(3727,66)(3801,73)(3870,79)
	(3933,86)(3991,93)(4044,101)
	(4091,108)(4134,115)(4172,122)
	(4207,129)(4237,136)(4264,143)
	(4288,150)(4310,157)(4329,164)
	(4348,170)(4365,177)(4381,184)
	(4398,190)(4414,197)(4432,204)
	(4451,211)(4471,218)(4494,225)
	(4519,232)(4546,239)(4578,247)
	(4612,254)(4651,262)(4694,269)
	(4741,277)(4792,285)(4848,292)
	(4908,300)(4971,307)(5038,315)
	(5108,322)(5180,328)(5253,335)
	(5325,340)(5421,346)(5513,351)
	(5598,355)(5677,357)(5750,358)
	(5817,358)(5878,357)(5934,355)
	(5986,352)(6033,349)(6077,345)
	(6118,340)(6155,336)(6190,330)
	(6223,325)(6253,319)(6281,313)
	(6307,307)(6331,301)(6353,295)
	(6372,290)(6389,285)(6404,280)
	(6417,276)(6427,273)(6435,270)
	(6441,268)(6445,267)(6448,266)
	(6449,265)(6450,265)
\path(525,1090)(529,1089)(537,1087)
	(550,1084)(569,1080)(593,1074)
	(621,1068)(649,1060)(678,1053)
	(705,1046)(729,1038)(751,1031)
	(771,1025)(788,1018)(803,1011)
	(816,1005)(827,997)(838,990)
	(846,983)(853,976)(860,968)
	(866,959)(871,951)(875,942)
	(878,932)(881,922)(882,913)
	(883,902)(882,892)(881,883)
	(878,873)(875,863)(871,854)
	(866,846)(860,837)(853,829)
	(846,822)(838,815)(827,808)
	(816,800)(803,794)(788,787)
	(771,780)(751,774)(729,767)
	(705,759)(678,752)(649,745)
	(621,737)(593,731)(569,725)
	(550,721)(537,718)(529,716)(525,715)
\thinlines
\path(6450,1390)(6449,1390)(6445,1387)
\path	(6389,1353)(6373,1342)(6359,1331)
\path	(6325,1290)(6320,1278)(6315,1265)
\path	(6307,1202)(6307,1186)(6309,1170)
\path	(6325,1115)(6331,1104)(6339,1094)
\path	(6389,1052)(6406,1041)(6423,1031)
\Thicklines
\thinlines
\path(3300,3790)(3301,3790)(3305,3787)
	(3314,3782)(3327,3774)(3344,3764)
	(3361,3753)(3377,3742)(3391,3731)
	(3402,3721)(3411,3711)(3419,3701)
	(3425,3690)(3430,3678)(3435,3665)
	(3438,3651)(3441,3635)(3443,3619)
	(3443,3602)(3443,3586)(3441,3570)
	(3438,3554)(3435,3540)(3430,3527)
	(3425,3515)(3419,3504)(3411,3494)
	(3402,3484)(3391,3474)(3377,3463)
	(3361,3452)(3344,3441)(3327,3431)
	(3314,3423)(3305,3418)(3301,3415)(3300,3415)
\thinlines
\path(3300,3790)(3299,3790)(3295,3787)
\path	(3239,3753)(3223,3742)(3209,3731)
\path	(3175,3690)(3170,3678)(3165,3665)
\path	(3157,3602)(3157,3586)(3159,3570)
\path	(3175,3515)(3181,3504)(3189,3494)
\path	(3239,3452)(3256,3441)(3273,3431)
\path(3300,3115)(3299,3115)(3295,3112)
\path	(3239,3078)(3223,3067)(3209,3056)
\path	(3175,3015)(3170,3003)(3165,2990)
\path	(3157,2927)(3157,2911)(3159,2895)
\path	(3175,2840)(3181,2829)(3189,2819)
\path	(3239,2777)(3256,2766)(3273,2756)
\Thicklines
\path(2766,1327)(2766,1326)(2767,1323)
	(2769,1318)(2771,1310)(2775,1301)
	(2780,1290)(2786,1278)(2793,1266)
	(2802,1253)(2813,1239)(2828,1225)
	(2845,1211)(2867,1196)(2894,1180)
	(2925,1165)(2955,1152)(2985,1142)
	(3011,1133)(3034,1125)(3052,1120)
	(3066,1116)(3077,1112)(3086,1110)
	(3095,1108)(3103,1106)(3112,1104)
	(3125,1102)(3141,1100)(3162,1098)
	(3189,1095)(3221,1093)(3259,1091)
	(3300,1090)(3341,1091)(3379,1093)
	(3411,1095)(3438,1098)(3460,1100)
	(3476,1102)(3488,1104)(3498,1106)
	(3507,1108)(3515,1110)(3525,1112)
	(3536,1116)(3550,1120)(3568,1125)
	(3590,1133)(3617,1142)(3646,1152)
	(3675,1165)(3705,1180)(3730,1196)
	(3751,1211)(3767,1225)(3779,1239)
	(3789,1253)(3797,1266)(3803,1278)
	(3807,1290)(3810,1301)(3813,1310)
	(3814,1318)(3815,1323)(3816,1326)(3816,1327)
\Thicklines
\thinlines
\path(3300,1090)(3299,1090)(3295,1087)
\path	(3239,1053)(3223,1042)(3209,1031)
\path	(3175,990)(3170,978)(3165,965)
\path	(3157,902)(3157,886)(3159,870)
\path	(3175,815)(3181,804)(3189,794)
\path	(3239,752)(3256,741)(3273,731)
\Thicklines
\path(2766,2503)(2766,2504)(2767,2507)
	(2769,2512)(2771,2520)(2775,2529)
	(2780,2540)(2786,2552)(2793,2564)
	(2802,2577)(2813,2591)(2828,2605)
	(2845,2619)(2867,2634)(2894,2650)
	(2925,2665)(2955,2678)(2985,2688)
	(3011,2697)(3034,2705)(3052,2710)
	(3066,2714)(3077,2718)(3086,2720)
	(3095,2722)(3103,2724)(3112,2726)
	(3125,2728)(3141,2730)(3162,2732)
	(3189,2735)(3221,2737)(3259,2739)
	(3300,2740)(3341,2739)(3379,2737)
	(3411,2735)(3438,2732)(3460,2730)
	(3476,2728)(3488,2726)(3498,2724)
	(3507,2722)(3515,2720)(3525,2718)
	(3536,2714)(3550,2710)(3568,2705)
	(3590,2697)(3617,2688)(3646,2678)
	(3675,2665)(3705,2650)(3730,2634)
	(3751,2619)(3767,2605)(3779,2591)
	(3789,2577)(3797,2564)(3803,2552)
	(3807,2540)(3810,2529)(3813,2520)
	(3814,2512)(3815,2507)(3816,2504)(3816,2503)
\path(525,3790)(526,3789)(528,3789)
	(531,3788)(534,3786)(540,3784)
	(546,3782)(555,3780)(565,3777)
	(577,3773)(592,3770)(608,3767)
	(626,3763)(648,3759)(671,3756)
	(698,3753)(728,3749)(763,3746)
	(801,3743)(845,3741)(895,3738)
	(952,3736)(1015,3734)(1087,3733)
	(1168,3732)(1256,3732)(1352,3732)
	(1416,3733)(1481,3733)(1546,3734)
	(1609,3736)(1670,3737)(1729,3739)
	(1784,3741)(1836,3743)(1884,3745)
	(1928,3748)(1968,3750)(2005,3753)
	(2038,3756)(2068,3759)(2095,3762)
	(2120,3765)(2142,3768)(2163,3771)
	(2181,3774)(2199,3777)(2216,3780)
	(2233,3783)(2251,3787)(2268,3790)
	(2287,3793)(2307,3796)(2329,3798)
	(2354,3801)(2381,3804)(2412,3806)
	(2446,3808)(2484,3811)(2526,3812)
	(2573,3814)(2625,3815)(2682,3816)
	(2744,3817)(2812,3817)(2884,3817)
	(2961,3816)(3042,3814)(3127,3813)
	(3214,3810)(3302,3807)(3390,3803)
	(3477,3799)(3562,3794)(3643,3789)
	(3719,3784)(3791,3779)(3857,3773)
	(3919,3767)(3975,3762)(4025,3756)
	(4071,3750)(4112,3744)(4148,3739)
	(4180,3733)(4208,3727)(4233,3722)
	(4256,3716)(4275,3711)(4293,3706)
	(4310,3700)(4325,3695)(4340,3690)
	(4354,3684)(4369,3679)(4385,3674)
	(4402,3668)(4421,3663)(4441,3658)
	(4465,3652)(4491,3647)(4520,3642)
	(4553,3636)(4590,3631)(4631,3625)
	(4677,3620)(4727,3615)(4781,3609)
	(4839,3604)(4902,3600)(4968,3595)
	(5037,3591)(5108,3588)(5180,3584)
	(5252,3582)(5343,3580)(5430,3579)
	(5512,3579)(5590,3580)(5661,3582)
	(5728,3585)(5790,3589)(5847,3593)
	(5901,3597)(5950,3602)(5997,3608)
	(6040,3614)(6081,3620)(6119,3627)
	(6155,3634)(6189,3641)(6221,3648)
	(6251,3655)(6279,3663)(6305,3670)
	(6329,3676)(6351,3683)(6371,3689)
	(6388,3694)(6404,3699)(6416,3703)
	(6427,3707)(6435,3710)(6441,3712)
	(6445,3713)(6448,3714)(6449,3715)(6450,3715)
\thinlines
\dottedline{90}(1800,1840)(1800,1240)
\thinlines
\path(6450,3715)(6449,3715)(6445,3712)
\path	(6389,3678)(6373,3667)(6359,3656)
\path	(6325,3615)(6320,3603)(6315,3590)
\path	(6307,3527)(6307,3511)(6309,3495)
\path	(6325,3440)(6331,3429)(6339,3419)
\path	(6389,3377)(6406,3366)(6423,3356)
\Thicklines
\path(6450,1015)(6446,1014)(6438,1012)
	(6425,1009)(6406,1005)(6382,999)
	(6354,993)(6326,985)(6297,978)
	(6270,971)(6246,963)(6224,956)
	(6204,950)(6187,943)(6172,936)
	(6159,930)(6148,922)(6138,915)
	(6129,908)(6122,901)(6115,893)
	(6109,884)(6104,876)(6100,867)
	(6097,857)(6094,847)(6093,838)
	(6092,827)(6093,817)(6094,808)
	(6097,798)(6100,788)(6104,779)
	(6109,771)(6115,762)(6122,754)
	(6129,747)(6138,740)(6148,733)
	(6159,725)(6172,719)(6187,712)
	(6204,705)(6224,699)(6246,692)
	(6270,684)(6297,677)(6326,670)
	(6354,662)(6382,656)(6406,650)
	(6425,646)(6438,643)(6446,641)(6450,640)
\Thicklines
\path(525,3790)(524,3790)(520,3787)
	(511,3782)(498,3774)(481,3764)
	(464,3753)(448,3742)(434,3731)
	(423,3721)(414,3711)(406,3701)
	(400,3690)(395,3678)(390,3665)
	(387,3651)(384,3635)(382,3619)
	(382,3602)(382,3586)(384,3570)
	(387,3554)(390,3540)(395,3527)
	(400,3515)(406,3504)(414,3494)
	(423,3484)(434,3474)(448,3463)
	(464,3452)(481,3441)(498,3431)
	(511,3423)(520,3418)(524,3415)(525,3415)
\path(525,3790)(526,3790)(530,3787)
	(539,3782)(552,3774)(569,3764)
	(586,3753)(602,3742)(616,3731)
	(627,3721)(636,3711)(644,3701)
	(650,3690)(655,3678)(660,3665)
	(663,3651)(666,3635)(668,3619)
	(668,3602)(668,3586)(666,3570)
	(663,3554)(660,3540)(655,3527)
	(650,3515)(644,3504)(636,3494)
	(627,3484)(616,3474)(602,3463)
	(586,3452)(569,3441)(552,3431)
	(539,3423)(530,3418)(526,3415)(525,3415)
\thinlines
\path(6450,2965)(6449,2965)(6445,2962)
\path	(6389,2928)(6373,2917)(6359,2906)
\path	(6325,2865)(6320,2853)(6315,2840)
\path	(6307,2777)(6307,2761)(6309,2745)
\path	(6325,2690)(6331,2679)(6339,2669)
\path	(6389,2627)(6406,2616)(6423,2606)
\Thicklines
\path(525,1840)(529,1839)(537,1837)
	(550,1834)(569,1830)(593,1824)
	(621,1818)(649,1810)(678,1803)
	(705,1796)(729,1788)(751,1781)
	(771,1775)(788,1768)(803,1761)
	(816,1755)(827,1747)(838,1740)
	(846,1733)(853,1726)(860,1718)
	(866,1709)(871,1701)(875,1692)
	(878,1682)(881,1672)(882,1663)
	(883,1652)(882,1642)(881,1633)
	(878,1623)(875,1613)(871,1604)
	(866,1596)(860,1587)(853,1579)
	(846,1572)(838,1565)(827,1558)
	(816,1550)(803,1544)(788,1537)
	(771,1530)(751,1524)(729,1517)
	(705,1509)(678,1502)(649,1495)
	(621,1487)(593,1481)(569,1475)
	(550,1471)(537,1468)(529,1466)(525,1465)
\Thicklines
\path(525,715)(526,715)(530,712)
	(539,707)(552,699)(569,689)
	(586,678)(602,667)(616,656)
	(627,646)(636,636)(644,626)
	(650,615)(655,603)(660,590)
	(663,576)(666,560)(668,544)
	(668,527)(668,511)(666,495)
	(663,479)(660,465)(655,452)
	(650,440)(644,429)(636,419)
	(627,409)(616,399)(602,388)
	(586,377)(569,366)(552,356)
	(539,348)(530,343)(526,340)(525,340)
\path(525,715)(526,715)(530,712)
	(539,707)(552,699)(569,689)
	(586,678)(602,667)(616,656)
	(627,646)(636,636)(644,626)
	(650,615)(655,603)(660,590)
	(663,576)(666,560)(668,544)
	(668,527)(668,511)(666,495)
	(663,479)(660,465)(655,452)
	(650,440)(644,429)(636,419)
	(627,409)(616,399)(602,388)
	(586,377)(569,366)(552,356)
	(539,348)(530,343)(526,340)(525,340)
\thinlines
\path(3300,415)(3301,415)(3305,412)
	(3314,407)(3327,399)(3344,389)
	(3361,378)(3377,367)(3391,356)
	(3402,346)(3411,336)(3419,326)
	(3425,315)(3430,303)(3435,290)
	(3438,276)(3441,260)(3443,244)
	(3443,227)(3443,211)(3441,195)
	(3438,179)(3435,165)(3430,152)
	(3425,140)(3419,129)(3411,119)
	(3402,109)(3391,99)(3377,88)
	(3361,77)(3344,66)(3327,56)
	(3314,48)(3305,43)(3301,40)(3300,40)
\path(3300,1090)(3301,1090)(3305,1087)
	(3314,1082)(3327,1074)(3344,1064)
	(3361,1053)(3377,1042)(3391,1031)
	(3402,1021)(3411,1011)(3419,1001)
	(3425,990)(3430,978)(3435,965)
	(3438,951)(3441,935)(3443,919)
	(3443,902)(3443,886)(3441,870)
	(3438,854)(3435,840)(3430,827)
	(3425,815)(3419,804)(3411,794)
	(3402,784)(3391,774)(3377,763)
	(3361,752)(3344,741)(3327,731)
	(3314,723)(3305,718)(3301,715)(3300,715)
\path(3300,3115)(3301,3115)(3305,3112)
	(3314,3107)(3327,3099)(3344,3089)
	(3361,3078)(3377,3067)(3391,3056)
	(3402,3046)(3411,3036)(3419,3026)
	(3425,3015)(3430,3003)(3435,2990)
	(3438,2976)(3441,2960)(3443,2944)
	(3443,2927)(3443,2911)(3441,2895)
	(3438,2879)(3435,2865)(3430,2852)
	(3425,2840)(3419,2829)(3411,2819)
	(3402,2809)(3391,2799)(3377,2788)
	(3361,2777)(3344,2766)(3327,2756)
	(3314,2748)(3305,2743)(3301,2740)(3300,2740)
\thicklines
\path(6450,3715)(6451,3715)(6455,3712)
	(6464,3707)(6477,3699)(6494,3689)
	(6511,3678)(6527,3667)(6541,3656)
	(6552,3646)(6561,3636)(6569,3626)
	(6575,3615)(6580,3603)(6585,3590)
	(6588,3576)(6591,3560)(6593,3544)
	(6593,3527)(6593,3511)(6591,3495)
	(6588,3479)(6585,3465)(6580,3452)
	(6575,3440)(6569,3429)(6561,3419)
	(6552,3409)(6541,3399)(6527,3388)
	(6511,3377)(6494,3366)(6477,3356)
	(6464,3348)(6455,3343)(6451,3340)(6450,3340)
\path(6450,2965)(6451,2965)(6455,2962)
	(6464,2957)(6477,2949)(6494,2939)
	(6511,2928)(6527,2917)(6541,2906)
	(6552,2896)(6561,2886)(6569,2876)
	(6575,2865)(6580,2853)(6585,2840)
	(6588,2826)(6591,2810)(6593,2794)
	(6593,2777)(6593,2761)(6591,2745)
	(6588,2729)(6585,2715)(6580,2702)
	(6575,2690)(6569,2679)(6561,2669)
	(6552,2659)(6541,2649)(6527,2638)
	(6511,2627)(6494,2616)(6477,2606)
	(6464,2598)(6455,2593)(6451,2590)(6450,2590)
\path(6450,1390)(6451,1390)(6455,1387)
	(6464,1382)(6477,1374)(6494,1364)
	(6511,1353)(6527,1342)(6541,1331)
	(6552,1321)(6561,1311)(6569,1301)
	(6575,1290)(6580,1278)(6585,1265)
	(6588,1251)(6591,1235)(6593,1219)
	(6593,1202)(6593,1186)(6591,1170)
	(6588,1154)(6585,1140)(6580,1127)
	(6575,1115)(6569,1104)(6561,1094)
	(6552,1084)(6541,1074)(6527,1063)
	(6511,1052)(6494,1041)(6477,1031)
	(6464,1023)(6455,1018)(6451,1015)(6450,1015)
\path(6450,640)(6451,640)(6455,637)
	(6464,632)(6477,624)(6494,614)
	(6511,603)(6527,592)(6541,581)
	(6552,571)(6561,561)(6569,551)
	(6575,540)(6580,528)(6585,515)
	(6588,501)(6591,485)(6593,469)
	(6593,452)(6593,436)(6591,420)
	(6588,404)(6585,390)(6580,377)
	(6575,365)(6569,354)(6561,344)
	(6552,334)(6541,324)(6527,313)
	(6511,302)(6494,291)(6477,281)
	(6464,273)(6455,268)(6451,265)(6450,265)
\Thicklines
\thinlines
\path(3300,415)(3299,415)(3295,412)
\path	(3239,378)(3223,367)(3209,356)
\path	(3175,315)(3170,303)(3165,290)
\path	(3157,227)(3157,211)(3159,195)
\path	(3175,140)(3181,129)(3189,119)
\path	(3239,77)(3256,66)(3273,56)
\path(6450,640)(6449,640)(6445,637)
\path	(6389,603)(6373,592)(6359,581)
\path	(6325,540)(6320,528)(6315,515)
\path	(6307,452)(6307,436)(6309,420)
\path	(6325,365)(6331,354)(6339,344)
\path	(6389,302)(6406,291)(6423,281)
\Thicklines
\path(525,3040)(524,3040)(520,3037)
 	(511,3032)(498,3024)(481,3014)
	(464,3003)(448,2992)(434,2981)
	(423,2971)(414,2961)(406,2951)
	(400,2940)(395,2928)(390,2915)
	(387,2901)(384,2885)(382,2869)
	(382,2852)(382,2836)(384,2820)
	(387,2804)(390,2790)(395,2777)
	(400,2765)(406,2754)(414,2744)
	(423,2734)(434,2724)(448,2713)
	(464,2702)(481,2691)(498,2681)
	(511,2673)(520,2668)(524,2665)(525,2665)
\path(525,1465)(524,1465)(520,1462)
 	(511,1457)(498,1449)(481,1439)
	(464,1428)(448,1417)(434,1406)
	(423,1396)(414,1386)(406,1376)
	(400,1365)(395,1353)(390,1340)
	(387,1326)(384,1310)(382,1294)
	(382,1277)(382,1261)(384,1245)
	(387,1229)(390,1215)(395,1202)
	(400,1190)(406,1179)(414,1169)
	(423,1159)(434,1149)(448,1138)
	(464,1127)(481,1116)(498,1106)
	(511,1098)(520,1093)(524,1090)(525,1090)
\path(525,715)(524,715)(520,712)
 	(511,707)(498,699)(481,689)
	(464,678)(448,667)(434,656)
	(423,646)(414,636)(406,626)
	(400,615)(395,603)(390,590)
	(387,576)(384,560)(382,544)
	(382,527)(382,511)(384,495)
	(387,479)(390,465)(395,452)
	(400,440)(406,429)(414,419)
	(423,409)(434,399)(448,388)
	(464,377)(481,366)(498,356)
	(511,348)(520,343)(524,340)(525,340)
\Thicklines
\path(525,2665)(529,2664)(537,2662)
 	(550,2659)(569,2655)(593,2649)
	(621,2643)(649,2635)(678,2628)
	(705,2621)(729,2613)(751,2606)
	(771,2600)(788,2593)(803,2586)
	(816,2580)(827,2572)(838,2565)
	(846,2558)(853,2551)(860,2543)
	(866,2534)(871,2526)(875,2517)
	(878,2507)(881,2497)(882,2488)
	(883,2477)(882,2467)(881,2458)
	(878,2448)(875,2438)(871,2429)
	(866,2421)(860,2412)(853,2404)
	(846,2397)(838,2390)(827,2383)
	(816,2375)(803,2369)(788,2362)
	(771,2355)(751,2349)(729,2342)
	(705,2334)(678,2327)(649,2320)
	(621,2312)(593,2306)(569,2300)
	(550,2296)(537,2293)(529,2291)(525,2290)
\path(525,3415)(529,3414)(537,3412)
 	(550,3409)(569,3405)(593,3399)
	(621,3393)(649,3385)(678,3378)
	(705,3371)(729,3363)(751,3356)
	(771,3350)(788,3343)(803,3336)
	(816,3330)(827,3322)(838,3315)
	(846,3308)(853,3301)(860,3293)
	(866,3284)(871,3276)(875,3267)
	(878,3257)(881,3247)(882,3238)
	(883,3227)(882,3217)(881,3208)
	(878,3198)(875,3188)(871,3179)
	(866,3171)(860,3162)(853,3154)
	(846,3147)(838,3140)(827,3133)
	(816,3125)(803,3119)(788,3112)
	(771,3105)(751,3099)(729,3092)
	(705,3084)(678,3077)(649,3070)
	(621,3062)(593,3056)(569,3050)
	(550,3046)(537,3043)(529,3041)(525,3040)
\path(6450,1765)(6446,1764)(6438,1762)
 	(6425,1759)(6406,1755)(6382,1749)
	(6354,1743)(6326,1735)(6297,1728)
	(6270,1721)(6246,1713)(6224,1706)
	(6204,1700)(6187,1693)(6172,1686)
	(6159,1680)(6148,1672)(6138,1665)
	(6129,1658)(6122,1651)(6115,1643)
	(6109,1634)(6104,1626)(6100,1617)
	(6097,1607)(6094,1597)(6093,1588)
	(6092,1577)(6093,1567)(6094,1558)
	(6097,1548)(6100,1538)(6104,1529)
	(6109,1521)(6115,1512)(6122,1504)
	(6129,1497)(6138,1490)(6148,1483)
	(6159,1475)(6172,1469)(6187,1462)
	(6204,1455)(6224,1449)(6246,1442)
	(6270,1434)(6297,1427)(6326,1420)
	(6354,1412)(6382,1406)(6406,1400)
	(6425,1396)(6438,1393)(6446,1391)(6450,1390)
\path(6450,2590)(6446,2589)(6438,2587)
 	(6425,2584)(6406,2580)(6382,2574)
	(6354,2568)(6326,2560)(6297,2553)
	(6270,2546)(6246,2538)(6224,2531)
	(6204,2525)(6187,2518)(6172,2511)
	(6159,2505)(6148,2497)(6138,2490)
	(6129,2483)(6122,2476)(6115,2468)
	(6109,2459)(6104,2451)(6100,2442)
	(6097,2432)(6094,2422)(6093,2413)
	(6092,2402)(6093,2392)(6094,2383)
	(6097,2373)(6100,2363)(6104,2354)
	(6109,2346)(6115,2337)(6122,2329)
	(6129,2322)(6138,2315)(6148,2308)
	(6159,2300)(6172,2294)(6187,2287)
	(6204,2280)(6224,2274)(6246,2267)
	(6270,2259)(6297,2252)(6326,2245)
	(6354,2237)(6382,2231)(6406,2225)
	(6425,2221)(6438,2218)(6446,2216)(6450,2215)
\path(6450,3340)(6446,3339)(6438,3337)
 	(6425,3334)(6406,3330)(6382,3324)
	(6354,3318)(6326,3310)(6297,3303)
	(6270,3296)(6246,3288)(6224,3281)
	(6204,3275)(6187,3268)(6172,3261)
	(6159,3255)(6148,3247)(6138,3240)
	(6129,3233)(6122,3226)(6115,3218)
	(6109,3209)(6104,3201)(6100,3192)
	(6097,3182)(6094,3172)(6093,3163)
	(6092,3152)(6093,3142)(6094,3133)
	(6097,3123)(6100,3113)(6104,3104)
	(6109,3096)(6115,3087)(6122,3079)
	(6129,3072)(6138,3065)(6148,3058)
	(6159,3050)(6172,3044)(6187,3037)
	(6204,3030)(6224,3024)(6246,3017)
	(6270,3009)(6297,3002)(6326,2995)
	(6354,2987)(6382,2981)(6406,2975)
	(6425,2971)(6438,2968)(6446,2966)(6450,2965)
\put(0,3565){\makebox(0,0)[lb]{\smash{{{\SetFigFont{12}{14.4}{\rmdefault}{\mddefault}{\updefault}$\alpha_1$}}}}}
\put(6750,3565){\makebox(0,0)[lb]{\smash{{{\SetFigFont{12}{14.4}{\rmdefault}{\mddefault}{\updefault}$\gamma_1$}}}}}
\put(0,490){\makebox(0,0)[lb]{\smash{{{\SetFigFont{12}{14.4}{\rmdefault}{\mddefault}{\updefault}$\alpha_r$}}}}}
\put(6675,340){\makebox(0,0)[lb]{\smash{{{\SetFigFont{12}{14.4}{\rmdefault}{\mddefault}{\updefault}$\gamma_t$}}}}}
\put(2775,190){\makebox(0,0)[lb]{\smash{{{\SetFigFont{12}{14.4}{\rmdefault}{\mddefault}{\updefault}$\beta_s$}}}}}
\put(2025,1465){\makebox(0,0)[lb]{\smash{{{\SetFigFont{12}{14.4}{\rmdefault}{\mddefault}{\updefault}$g_1$}}}}}
\put(2775,3565){\makebox(0,0)[lb]{\smash{{{\SetFigFont{12}{14.4}{\rmdefault}{\mddefault}{\updefault}$\beta_1$}}}}}
\put(5400,1840){\makebox(0,0)[lb]{\smash{{{\SetFigFont{12}{14.4}{\rmdefault}{\mddefault}{\updefault}$g_2$}}}}}
\Thicklines
\path(525,1465)(526,1465)(530,1462)
	(539,1457)(552,1449)(569,1439)
	(586,1428)(602,1417)(616,1406)
	(627,1396)(636,1386)(644,1376)
	(650,1365)(655,1353)(660,1340)
	(663,1326)(666,1310)(668,1294)
	(668,1277)(668,1261)(666,1245)
	(663,1229)(660,1215)(655,1202)
	(650,1190)(644,1179)(636,1169)
	(627,1159)(616,1149)(602,1138)
	(586,1127)(569,1116)(552,1106)
	(539,1098)(530,1093)(526,1090)(525,1090)
\end{picture}
}

%% file: TQFTpost.bbl
\begin{thebibliography}{10}

\bibitem{Abrams}
Lowell Abrams.
\newblock Two-dimensional topological quantum field theories and {F}robenius
  algebras.
\newblock {\em J. Knot Theory Ramifications}, 5(5):569--587, 1996.

\bibitem{Bryan-Pandharipande-in-prep}
Jim Bryan and Rahul Pandharipande.
\newblock In preparation.

\bibitem{Br-Pa-rigidity}
Jim Bryan and Rahul Pandharipande.
\newblock {Rigidity of curves in Calabi-Yau 3-folds}.
\newblock In preparation.

\bibitem{Br-Pa}
Jim Bryan and Rahul Pandharipande.
\newblock B{P}{S} states of curves in {C}alabi-{Y}au 3-folds.
\newblock {\em Geom. Topol.}, 5:287--318 (electronic), 2001.
\newblock preprint version: math.AG/0009025.

\bibitem{Dijkgraaf-mirror95}
Robbert Dijkgraaf.
\newblock Mirror symmetry and elliptic curves.
\newblock In {\em The moduli space of curves (Texel Island, 1994)}, volume 129
  of {\em Progr. Math.}, pages 149--163. Birkh\"auser Boston, Boston, MA, 1995.

\bibitem{Dijkgraaf-Witten90}
Robbert Dijkgraaf and Edward Witten.
\newblock Topological gauge theories and group cohomology.
\newblock {\em Comm. Math. Phys.}, 129(2):393--429, 1990.

\bibitem{EGH}
Y.~Eliashberg, A.~Givental, and H.~Hofer.
\newblock Introduction to symplectic field theory.
\newblock {\em Geom. Funct. Anal.}, (Special Volume, Part II):560--673, 2000.
\newblock GAFA 2000 (Tel Aviv, 1999).

\bibitem{Faber-Pandharipande-03}
C.~Faber and R.~Pandharipande.
\newblock {Relative maps and tautological classes}.
\newblock arXiv:math.AG/0304485.

\bibitem{Fa-Pa}
C.~Faber and R.~Pandharipande.
\newblock Hodge integrals and {G}romov-{W}itten theory.
\newblock {\em Invent. Math.}, 139(1):173--199, 2000.

\bibitem{Fantechi-Pandharipande}
B.~Fantechi and R.~Pandharipande.
\newblock Stable maps and branch divisors.
\newblock {\em Compositio Math.}, 130(3):345--364, 2002.
\newblock Preprint: math.AG/9905104.

\bibitem{Freed-Quinn}
Daniel~S. Freed and Frank Quinn.
\newblock Chern-{S}imons theory with finite gauge group.
\newblock {\em Comm. Math. Phys.}, 156(3):435--472, 1993.

\bibitem{Graber-Vakil}
Tom Graber and Ravi Vakil.
\newblock {Relative virtual localization and vanishing of tautological classes
  on moduli spaces of curves}.
\newblock arXiv:math.AG/0309227.

\bibitem{Ionel-Parker00}
Eleny-Nicoleta Ionel and Thomas~H. Parker.
\newblock {The Symplectic Sum Formula for Gromov-Witten Invariants}.
\newblock arXiv:math.SG/0010217.

\bibitem{Ionel-Parker-Annals2003}
Eleny-Nicoleta Ionel and Thomas~H. Parker.
\newblock Relative {G}romov-{W}itten invariants.
\newblock {\em Ann. of Math. (2)}, 157(1):45--96, 2003.

\bibitem{Katz-Liu}
Sheldon Katz and Chiu-Chu~Melissa Liu.
\newblock Enumerative geometry of stable maps with {L}agrangian boundary
  conditions and multiple covers of the disc.
\newblock {\em Adv. Theor. Math. Phys.}, 5(1):1--49, 2001.

\bibitem{Kock:FA-2DTQFT}
Joachim Kock.
\newblock Frobenius algebras and {2D} topological quantum field theories.
\newblock To appear in LMSST, Cambridge University Press.

\bibitem{Li-Ruan}
An-Min Li and Yongbin Ruan.
\newblock Symplectic surgery and {G}romov-{W}itten invariants of {C}alabi-{Y}au
  3-folds.
\newblock {\em Invent. Math.}, 145(1):151--218, 2001.

\bibitem{Li-relative1}
Jun Li.
\newblock Stable morphisms to singular schemes and relative stable morphisms.
\newblock {\em J. Differential Geom.}, 57(3):509--578, 2001.
\newblock arXiv:math.AG/0009097.

\bibitem{Li-relative2}
Jun Li.
\newblock A degeneration formula of {GW}-invariants.
\newblock {\em J. Differential Geom.}, 60(2):199--293, 2002.

\bibitem{Li-Song}
Jun Li and Yun~S. Song.
\newblock Open string instantons and relative stable morphisms.
\newblock {\em Adv. Theor. Math. Phys.}, 5(1):67--91, 2001.
\newblock arXiv:hep-th/0103100.

\bibitem{Moore-ITPlecture}
Greg Moore.
\newblock {ITP lectures on branes, K-theory, and RR-charges}, 2001.
\newblock {Available online at }\\{
  http://online.kitp.ucsb.edu/online/mp01/moore1/.}

\bibitem{Okounkov-Pandharipande-completed-cycles}
Andrei Okounkov and Rahul Pandharipande.
\newblock {Gromov-Witten theory, Hurwitz theory, and completed cycles}.
\newblock arXiv:math.AG/0204305.

\bibitem{Pandharipande-degenerate-contributions}
R.~Pandharipande.
\newblock Hodge integrals and degenerate contributions.
\newblock {\em Comm. Math. Phys.}, 208(2):489--506, 1999.

\bibitem{Quinn}
Frank Quinn.
\newblock Lectures on axiomatic topological quantum field theory.
\newblock In {\em Geometry and quantum field theory (Park City, UT, 1991)},
  volume~1 of {\em IAS/Park City Math. Ser.}, pages 323--453. Amer. Math. Soc.,
  Providence, RI, 1995.

\bibitem{Sawin}
Stephen Sawin.
\newblock Direct sum decompositions and indecomposable {TQFT}s.
\newblock {\em J. Math. Phys.}, 36(12):6673--6680, 1995.

\end{thebibliography}
